\documentclass[times, doublespace]{fldauth}
\usepackage{graphicx}
\usepackage{amssymb,amsxtra,amsmath}
\usepackage{latexsym}
\usepackage{epstopdf}

\usepackage{verbatim}
\usepackage{enumerate}
\usepackage{verbatim}
\usepackage{ifthen}
\usepackage{calc,fancyhdr}
\usepackage{indentfirst}
\usepackage{parskip}
\usepackage{cite}

\usepackage{color}
\definecolor{black}{rgb}{0,0,0}

\definecolor{none}{rgb}{1,1,1}

\definecolor{red}{rgb}{1.,.0,.0}

\definecolor{blue}{rgb}{0,0,1}

\renewcommand{\vec}[1]{\underline{#1}}

\newcommand{\ignore}[1]{}

\long\def\comment#1{}

\def \vu{ \boldsymbol{u} }
\def \vU{ \boldsymbol{U} }

\def \vx{ \boldsymbol{x} }

\def \vf{ \boldsymbol{f} }
\def \vX{ \boldsymbol{X} }
\def \vZ{ \boldsymbol{Z} }
\def \vXd{ \boldsymbol{X}^{\rm d} }

\def \vF{ \boldsymbol{F} }

\def \Xd{ X^{\rm d} }

\def \Yd{ Y^{\rm d} }

\def \ucd{ \vec{c}^X_{\rm d} }
\def \uXd{ \vec{X}_{\rm d}(t) }
\def \uvXd{ \vec{\vX}_{\rm d} }
\def \uvXs{ \vec{\vX}_{\rm s} }
\def \uYd{ \vec{Y}_{\rm d}(t) }
\def \uXs{ \vec{X}_{\rm s}(t) }
\def \uTd{ \vec{T}_{\rm d} }
\def \uvZd{ \vec{\vZ}_{\rm d} }
\def \uvFs{ \vec{\vF}_{\rm s} }
\def \uvFTs{ \vec{\vF}_{\rm s}^{\rm T} }
\def \uvFBs{ \vec{\vF}_{\rm s}^{\rm B} }
\def \uvFCs{ \vec{\vF}_{\rm s}^{\rm C} }
\def \tngtd{ \tngt_{\rm d} }
\def \utngtd{ \vec{\tngt}_{\rm d} }
\def \uutngtd{ \vec{\hat{\tngt}}_{\rm d} }

\def \Es{ \mathcal{E}_{\rm s} }
\def \ub{ \vec{b}_{\lon} }
\def \ubd{ \ub^{n} }

\newcommand{\Dlond}[1]{\mathcal{D}_{\lond}^{#1}}
\newcommand{\Dlons}[1]{\mathcal{D}_{\lons}^{#1}}

\newcommand{\lon}{\lambda}

\newcommand{\lond}{\lon^{\rm d}}

\newcommand{\lons}{\lambda^{\rm s}}

\newcommand{\Nd}{N_{\rm d}}
\newcommand{\Ns}{N_{\rm s}}
\newcommand{\lf}{\left}
\newcommand{\rt}{\right}
\newcommand{\ep}{\varepsilon}
\newcommand{\tngt}{\boldsymbol{\tau}}
\newcommand{\utngt}{\hat{\boldsymbol{\tau}}}

\newcommand{\frc}{\vF}
\newcommand{\frcT}{\frc^{\rm T}}
\newcommand{\frcB}{\frc^{\rm B}}
\newcommand{\frcC}{\frc^{\rm C}}
\newcommand{\kt}{k_{\rm t}}
\newcommand{\kb}{k_{\rm b}}
\newcommand{\KC}{K_{\rm C}}
\newcommand{\ellC}{l_{0,\rm C}}

\begin{document}

\title{Augmenting the Immersed Boundary Method with Radial Basis Functions (RBFs) for the Modeling of Platelets in Hemodynamic Flows}

\author{Varun Shankar\affil{1}\corrauth, Grady B. Wright\affil{2}, Robert M. Kirby\affil{3}, and Aaron L. Fogelson\affil{4}}

\address{\affilnum{1}Department of Mathematics, Univ.\ of Utah, Salt Lake City, UT, USA\break
\affilnum{2}Department of Mathematics, Boise State Univ., Boise, ID, USA
\break\affilnum{3}School of Computing, Univ.\ of Utah, Salt Lake City, UT, USA
\break\affilnum{4}Departments of Mathematics and Bioengineering, Univ.\ of Utah, Salt Lake City, UT, USA}

\corraddr{Department of Mathematics, Univ.\ of Utah, Salt Lake City, UT, USA. E-mail: vshankar@math.utah.edu}


\bigskip
\begin{abstract}
We present a new computational method by extending the Immersed Boundary (IB) method with a geometric model based on parametric Radial Basis Function (RBF) interpolation of the Lagrangian structures. Our specific motivation is the modeling of platelets in hemodynamic flows, though we anticipate that our method will be useful in other applications involving surface elasticity. The efficacy of our new RBF-IB method is shown through a series of numerical experiments. Specifically, we test the convergence of our method and compare our method with the traditional IB method in terms of computational cost, maximum stable time-step size and volume loss. We conclude that the RBF-IB method has advantages over the traditional Immersed Boundary method, and is well-suited for modeling of platelets in hemodynamic flows.
\end{abstract}

\maketitle

\section{Introduction}
\label{sec:intro}

The Immersed Boundary (IB) Method was introduced by Charles Peskin 
in the early 1970's to solve the coupled equations of motion of a
viscous, incompressible fluid and one or more massless, elastic
surfaces or objects immersed in the fluid~\cite{Peskin:1977}. The IB method was originally developed to model blood flow in the
heart and through heart valves~\cite{Peskin:1977,Peskin:1980,Peskin:1989}, 
but has since been used in a wide variety of other applications, particularly in
biofluid dynamics problems where complex geometries and immersed
elastic membranes or structures are present and make traditional
computational approaches difficult. Examples include swimming of organisms~\cite{Fauci_Fogelson:1993,Fauci_Peskin:1988}, biofilm
processes~\cite{Dillon:1996}, mechanical properties of
cells~\cite{Agresar:1998}, cochlear dynamics~\cite{Beyer:1992},
and insect flight~\cite{Miller_Peskin:2004b,Miller_Peskin:2004}. In this work, we are motivated by the application of the IB method to platelet aggregation in blood clotting, but expect
our method to be useful in other applications as well.

Intravascular blood clots (thrombi) are initiated by damage to the
endothelial cell lining of a blood vessel and involve the formation on
the damaged surface of clumps of cells intermixed with a fibrous
protein gel.  The cells involved in this process are platelets, and
the subject of this paper is a new approach to modeling platelets in
order to simulate their adhesion to the injured vascular wall and
cohesion with one another during the formation of a thrombus. The IB method is used to describe the mechanical interactions among a collection of discrete platelets, the background fluid, and the vessel wall. However, we now model our platelets with Radial basis functions (RBFs) in order to achieve more accurate and less costly simulations.

In this introduction, we briefly describe the relevant biology, describe
how the IB method has been used in our previous platelet aggregation
simulations, and give an overview of how use of the new method changes
this description.

\subsection{Modeling the mechanics of platelet aggregation}
\label{sec:micromodels}

Disruption of the endothelial cell lining exposes collagen and adsorbed von Willebrand factor (vWF) molecules in the subendothelial matrix to the blood.  Platelets adhere to both molecules via specific receptor molecules on the platelets' surfaces.  In addition to slowing or stopping platelet motion over the subendothelium, this
binding triggers intracellular signaling pathways that lead to
platelet activation \cite{JACKSON:2003:SEU,SAVAGE:1996:IPA}.  

Our platelet aggregation models \cite{Fauci_Fogelson:1993,Fogelson:1984,Fogelson2008,FOGELSON:2003:CMB,YU:2000:TMP} track the motion and behavior of a collection of individual platelets
as they interact with the suspending fluid, one another, and the
vessel walls.  We also track fluid concentrations of platelet
activating chemicals, cell-cell and cell-surface forces, fluid motion,
and the local fluid forces on the growing thrombus.  In our models,
nonactivated platelets are activated by proximity to reactive sites on
the injured wall, or through exposure to a sufficiently high
concentration of activator in the fluid.  Activation enables a
platelet to cohere with other activated platelets and to secrete
additional activator.  The platelets and the secreted chemical move by
advection with the fluid and diffusion relative to it.  Each platelet
is represented as an IB object, {\emph i.e.,} as a collection of
elastically-linked Lagrangian points that each move at the local fluid
velocity.  New elastic links are created dynamically to model the
adhesion of a platelet to the injured wall or the cohesion of
activated platelets to one another.  Multiple links can form between a
pair of activated model platelets or between a model platelet and the
injured wall, and these links collectively represent the ensemble of
molecular bridges binding real platelets to one another or to the
damaged vessel.  The links exert forces on the surrounding fluid to
resist motions which would otherwise separate the linked entities.
Links may break if subject to sufficiently high stress.  Model variables are fully coupled: the fluid carries the
activator and platelets, while the interplatelet forces, potentiated
by chemically-induced activation of the platelets, determine the local
flow.  In this paper, we focus on mechanical interactions, not the
activation process, and so we specify the conditions under which a
platelet becomes activated and able to cohere with other activated
platelets.

\subsection{Motivation for the RBF-IB method}
\label{subsec:introrbfibm}

We model platelets as closed curves of interconnected IB points in 2D.  A platelet's
area or volume is determined by the region enclosed by the curve or
surface and is preserved because of the incompressibility of the
fluid. Inactive platelets are approximately elliptical in 2D models, while activated platelets are approximately circular. Piecewise linear approximations of platelets are currently used in IB methods applied to the simulation of platelet
aggregation ({\em e.g.}\cite{Fauci_Fogelson:1993,Fogelson:1984,Fogelson2008}).

In previous work \cite{Shankar2012}, we found that interpolation with radial basis functions restricted to the circle (or sphere in 3D)
offered accuracy and computational cost comparable to that offered by Fourier-based methods in modeling an infinitely smooth target shape, its normals and tension forces computed on its surface. Furthermore, interpolation with radial basis functions resulted in better convergence (often an order more) than that offered by both Fourier-based methods when the target shape had only one or two underlying derivatives. In general, use of radial basis functions led to a computational cost comparable to that of Fourier-based methods
and orders of magnitude lower (for the same accuracy) than that of the standard combination of techniques (piecewise quadratic interpolation and finite differences) used in many IB methods. This RBF based geometric model has since been used in a variety of applications~\cite{SWFKIJNMF2014,SWFKJSC2014,FuselierWright2013}. It is currently being extended to the representation of open elastic curves immersed in a purely viscous fluid for use within the method of Regularized Stokeslets~\cite{SOIJNMF2015}.

We now turn our attention to exploring
the consequences of using this parametric RBF geometric model within the full IB method,
with platelet aggregation as our target application. We seek to determine if the advantages inherent in the RBF interpolation of static shapes carries over to full-fledged IB simulations, and also if the RBF interpolation can give us benefits that are apparent only in full-fledged IB simulations. In this work, we propose a new immersed boundary algorithm that utilizes the features afforded by our RBF geometric model.

The paper is organized as follows. In Section \ref{sec:IBreview} we briefly
discuss the traditional Immersed Boundary Method for simulating
fluid-structure interaction. In Section \ref{sec:GeometricModeling} we
review the piecewise linear and RBF geometric modeling strategies and review the
components necessary for handling immersed elastic structures in the
IB method. In Section \ref{sec:implementation} we provide
details of the spatial and temporal discretizations of both versions of the IB method. In Section
\ref{sec:results} we present our comparison of the RBF-IB method with the traditional
IB method in terms of convergence, accuracy, area loss and time-step size. We also provide energy estimates for RBF-IB simulations. We then present results from a large platelet aggregation simulation in 2D. Section \ref{sec:discussion} contains a summary of our findings and a discussion of future research directions.

\underline{Notation}: We denote vectors with as many components as the spatial dimension in bold. We denote vectors with as many components as the number of data sites ($\Nd$) or sample sites ($\Ns$) by underlining. We indicate matrices with ($\Nd$) or ($\Ns$) rows and two columns in bold with underlining.
\section{Review of the Immersed Boundary Method}
\label{sec:IBreview}

To review the IB method, we focus on a simple
two-dimensional model problem in which a single fluid-filled closed
elastic membrane is immersed in a viscous fluid. The physics of the model problem is that an elastic membrane is under
tension and exerts forces on the adjacent fluid. These forces may
cause the fluid to move and, correspondingly, cause points on the membrane to move
along with the fluid. In the IB method, the fluid is described in
the Eulerian frame through a velocity field $\vu(\vx,t)$ and pressure
field $p(\vx,t)$ defined at every point $\vx$ in the physical domain
$\Omega$.  The elastic membrane is described in the Lagrangian frame.  Let
the elastic membrane be parameterized by $q \in \Gamma$, and denote by $\vX(q,t)$
the spatial coordinates at time $t$ of the membrane point labeled by
$q$.  The IB equations are the following coupled equations of motion
for the fluid variables $\vu(\vx,t)$ and $p(\vx,t)$ and the membrane
configuration $\vX(q,t)$.
\begin{equation}
\rho (\vu_t + \vu \cdot \nabla \vu ) =  -\nabla p + \mu \nabla^2 \vu  +
\vf,  \hspace{.25in} \nabla \cdot \vu = 0,
\label{eq:nseq}
\end{equation}
\begin{equation}
\vF(q,t) ~=~ \boldsymbol{\mathit{F}}\lf(\vX (q,t), \frac{\partial}{\partial q}\vX(q,t)\rt),
\label{eq:ibF}
\end{equation}
\begin{equation}
\vf(\vx,t) ~=~ \int_{\Gamma} \vF (q,t) ~\delta( \vx - \vX(q,t) )~dq,
\label{eq:fluidf}
\end{equation}
\begin{equation}
\frac{\partial \vX}{\partial t}(q,t) = \int_{\Omega} \vu(\vx,t)
~\delta(\vx-\vX(q,t))d\vx.
\label{eq:ibU}
\end{equation}

Equations \eqref{eq:nseq} are the Navier Stokes equations which describe the
dynamics of a viscous incompressible fluid, of constant density $\rho$
and constant viscosity $\mu$, driven by a force density $\vf$ which
here arises because of the elastic deformation of the immersed
membrane.  Equation \eqref{eq:ibF} specifies the elastic force (per unit
$q$) at each point of the immersed boundary object.  The functional
dependence of this force on the state of the boundary is specified
appropriately to the material being modeled. Equation \eqref{eq:fluidf} defines the fluid force density
$\vf(\vx,t)$ in terms of the immersed boundary elastic force density
$\vF$. Equation \eqref{eq:ibU} specifies that the velocity of each immersed
boundary point equals the fluid velocity at the same location, a formulation of the no-slip boundary condition for viscous flows.
In the model problem and the platelet applications, we assume that the IB objects are neutrally buoyant; the
IB membrane itself carries no mass and each object's mass is attributed to the fluid in which it sits. For more on the IB method, see \cite{Peskin:2002}.
\section{Geometric modeling of platelets}
\label{sec:GeometricModeling}

In this section, we review the two geometric modeling strategies to be
compared in the context of the Immersed Boundary method applied to platelet
aggregation. For a full description
of these strategies, see \cite{Shankar2012}.

\subsection{Piecewise linear model}
\label{sec:Geom_PL}

In the traditional IB method, the surfaces of platelets are represented by a collection of Immersed Boundary points. We henceforth refer to the IB points in the traditional IB method as \emph{sample sites}, and denote them by $\vX_{\rm s}(t) = \vX(q,t)$ for  each discrete $q \in \Gamma$ and at a particular time $t$. The surface elastic forces of the platelets are spread from these sample sites into the neighboring fluid. Both tension and bending forces are computed using a finite difference discretization of force models at sample sites. An explicit piecewise linear interpolant of the surface is not formed. If other information (such as normal vectors) is needed at the sample points, an approximation to the surface represented by the sample sites may be formed by a piecewise quadratic interpolation of the sample sites ({\em e.g.},~\cite{YaoFogelson2012}). After the incompressible Navier Stokes equations are solved, velocities from the portions of the Eulerian grid
surrounding the sample sites are interpolated to the sample sites using a discretization of Equation \ref{eq:ibU} and used to move the platelets.
 
\subsection{Parametric RBF model}
\label{sec:Geom_RBF}

The RBF method is a popular tool for approximating multidimensional
scattered data. For an overview of the theory and application of
this method, see the books by Fasshauer~\cite{Fasshauer:2007} and Wendland~\cite{Wendland:2004}. The
restriction of the RBF method to interpolation on a circle and/or
sphere is discussed by Fasshauer and Schumaker~\cite[\S
6]{FasshauerSchumaker:1998}.  When restricted to
these domains, the RBF method is referred to as the \emph{spherical basis function}
(SBF) method~\cite[Ch. 17]{Wendland:2004}. Several studies have provided error
estimates for RBF interpolation on circles and spheres; in fact, these interpolants can provide
spectral accuracy provided the underlying target function is sufficiently smooth~\cite{JetterStocklerWard:1999,NarcSunWard:2007}. The RBF method has also been used successfully for numerically solving partial differential equations on the surface of the sphere~\cite{FlyerWright:2007,FlyerWright:2009}, as well as more general surfaces~\cite{FuselierWright2013,Piret2012}.

Here, we present the RBF model developed in our earlier work~\cite{Shankar2012}. 
It is based on explicit parametric representations of the
objects. Since our target objects are platelets, which in 2D models are nearly elliptical or
circular, we choose a polar parameterization. We use our model to define
operators necessary for the computation of geometric and mechanical quantities required by the IB method.

We represent a platelet surface at any time $t$ parametrically by
\begin{equation}
\label{eq:2D_obj}
\vX(\lon,t) = (X(\lon,t),Y(\lon,t))
\end{equation}
where $0 \leq \lon \leq 2\pi$ is the parametric variable and $\vX(0,t)=\vX(2\pi,t)$. We explicitly track
a finite set of $\Nd$ points $\vXd_1(t),\ldots,\vXd_{\Nd}(t)$, which we refer to as \emph{data sites}.  Here $\vXd_j(t) := \vX(\lond_j,t)$, $j=1,\ldots,\Nd$, and we refer to the parametric coordinates $\lond_1,\ldots,\lond_{\Nd}$ as the \emph{data site nodes} (or simply \emph{nodes}).   We construct each component of $\vX$ by using a smooth RBF interpolant of the data sites in parameter space as discussed in detail below. We also make use of derivatives of the interpolant at the data sites and we use the interpolant and its derivatives at another
set of prescribed \emph{sample points} or \emph{sample sites}, which correspond to $\Ns$ parameter values: $\lons_1,...,\lons_{\Ns}$. 

We first explain how to construct an RBF interpolant to the $X$ component of $\vX$ using the data $(\lond_1,\Xd_1(t)),...,(\lond_{\Nd},\Xd_{\Nd}(t))$; the construction of the $Y$ component follows in a similar manner.  Let $\phi(r)$ be
a scalar-valued radial kernel, whose choice we discuss below. Define $X(\lon,t)$ by
\begin{align}
X(\lon,t) &= \sum_{k=1}^{N_d} c^X_k \phi\lf(\sqrt{2 - 2\cos(\lon-\lond_k)}\rt).
\label{eq:circ_rbf_interp}
\end{align}
Note that the square root term in Equation \eqref{eq:circ_rbf_interp} is
the Euclidean distance between the points on the unit circle whose angular coordinates are
$\lon$ and $\lond_k$. We have found that the distance argument $r = \sqrt{2 - 2\cos(\lon-\lond_k)}$  is far more accurate
for the geometric of modeling static closed curves and surfaces than, say, $r = | \lon - \lond_k|$. In addition, recent work has shown that the periodic distance argument gives results identical to those given by the non-periodic distance argument for the modeling of both static \emph{and} dynamic open curves~\cite{SOIJNMF2015}. While other distance arguments could be considered (for example, using geodesic distance in place of Euclidean distance), Fuselier and Wright have shown that RBF interpolation can produce favorable error estimates in the interpolation of functions on submanifolds of $\mathbb{R}^n$ even when no knowledge of the surface metric is used~\cite{FuselierWright:2010}. For these reasons, we restrict our attention to the periodic distance argument in \eqref{eq:circ_rbf_interp}. For this paper, we use the multiquadric (MQ) radial kernel function, given by
\begin{align}
\text{MQ:}\quad & \phi(r) = \sqrt{1 + (\ep r)^2}, \label{eq:mq}
\end{align}
where $\ep$ is called the shape parameter. The choice of $\ep$ is discussed in Section 5. To have $X(\lon,t)$ interpolate the
given data, we require that the coefficients $c_k^X, k=1,...,\Nd$ satisfy the
following system of equations:
\begin{align}
\underbrace{
\begin{bmatrix}
\phi\lf(r_{1,1}\rt) & \cdots & \phi\lf(r_{1,\Nd}\rt) \\
\phi\lf(r_{2,1}\rt) & \cdots & \phi\lf(r_{2,\Nd}\rt) \\
\vdots & \ddots & \vdots \\
\phi\lf(r_{\Nd,1}\rt) & \cdots & \phi\lf(r_{\Nd,\Nd}\rt)
\end{bmatrix}}_{\displaystyle A}
\underbrace{
\begin{bmatrix}
c^X_1 \\
c^X_2 \\
\vdots \\
c^X_{\Nd}
\end{bmatrix}}_{\displaystyle \ucd}
= 
\underbrace{
\begin{bmatrix}
\Xd_1(t) \\
\Xd_2(t) \\
\vdots \\
\Xd_{\Nd}(t)
\end{bmatrix}}_{\displaystyle \uXd },
\label{eq:rbf_linsys}
\end{align}
where $r_{j,k} = \sqrt{2 - 2\cos(\lond_j - \lond_k)}$. Since $r_{j,k} = r_{k,j}$, the matrix
$A$ in this system is symmetric. More importantly, for the MQ kernels, $A$ is non-singular, with the global support and infinite smoothness of $\phi(r)$ lending itself to spectral accuracy and convergence on smooth problems~\cite{Fasshauer:2007, Wendland:2004}. One could alternatively use any of the other infinitely-smooth kernels like the Gaussian (GA) or the Inverse Multiquadric (IMQ) in place of the MQ kernel.

In our application, we want to be able to evaluate $X(\lon,t)$ at sample sites corresponding to parameter values $\lons_1,...,\lons_{\Ns}$, that stay fixed over time.  While we could use Equation \eqref{eq:circ_rbf_interp} to do this, it is much more convenient from a notational and computational perspective to construct an \emph{evaluation matrix} that combines the linear operations of constructing the interpolant to $\vec{\vX}_{\rm d}(t)=[\uXd,\; \uYd]$, for any $t$, and evaluating it at $\lons_1,...,\lons_{\Ns}$.  The evaluation matrix can be constructed by first noting that Equation \eqref{eq:circ_rbf_interp} can be written as
\begin{align}
X(\lon,t) = 
\underbrace{
\begin{bmatrix}
\phi\lf(\sqrt{2 - 2\cos(\lon-\lond_1})\rt) & \cdots  & \phi\lf(\sqrt{2 - 2\cos(\lon-\lond_{\Nd})}\rt)
\end{bmatrix}}_{\displaystyle \vec{b}(\lon)^T}
\ucd.
\label{eq:b_vec}
\end{align}
Since $\ucd = A^{-1} \uXd$, we can write Equation \eqref{eq:circ_rbf_interp} as $X(\lon,t) = \vec{b}(\lon)^{T}A^{-1} \uXd$.  The evaluation of $X(\lon,t)$ at $\lons_1,...,\lons_{\Ns}$ can then be obtained as follows:
\begin{align}
\underbrace{
\begin{bmatrix}
X(\lons_1,t) \\
\vdots \\
X(\lons_{\Ns},t)
\end{bmatrix}}_{\displaystyle \uXs}
= 
\underbrace{
\begin{bmatrix}
\vec{b}(\lons_1)^T \\
\vdots \\
\vec{b}(\lons_{\Ns})^T
\end{bmatrix}}_{\displaystyle B}
A^{-1}\uXd
= 
\underbrace{BA^{-1}}_{\displaystyle \Es} \uXd.
\label{eq:rbf_eval}
\end{align}
So, given the data sites $\uXd$ at any time $t$, we can interpolate their coordinates with an RBF expansion \emph{and} evaluate the interpolant at the sample site nodes $\lons_1,...,\lons_{\Ns}$ to get $\uXs$ by the matrix-vector product $\Es\uXd$.  In fact, this same procedure can be used to give values at sample site nodes for any quantity whose values we have at data site nodes and which
we represent using an RBF expansion ({\em e.g.}, $\uYd = [\Yd_1(t)  \cdots  \Yd_{\Nd}(t)]^{T}$).  Furthermore, the evaluation matrix $\Es$ can be precomputed once at $t=0$ and used for all subsequent times.

We also need to compute geometric quantities such as tangent vectors, and mechanical quantities such as forces at data sites and/or sample sites. These quantities require computing derivatives with respect to $\lon$ of the platelet surface coordinates $\lf(X(\lon,t),Y(\lon,t)\rt)$.  We use the RBF-based representation of the surface to compute these derivatives, and we will express derivatives of the RBF interpolant in matrix-vector form.   Toward this end, we use similar notation to Equation \eqref{eq:b_vec} and define the vector
\begin{align*}
\ubd(\tilde{\lon}) :=&  \lf.\dfrac{\partial^n}{\partial \lon^n}\vec{b}(\lon)\rt|_{\lon=\tilde{\lon}} \\
=&  
\begin{bmatrix}
\lf.\dfrac{\partial^n}{\partial \lon^n}\phi\lf(\sqrt{2 - 2\cos(\lon-\lond_1})\rt)\rt|_{\lon = \tilde{\lon}} & \cdots  & \lf.\dfrac{\partial^n}{\partial \lon^n}\phi\lf(\sqrt{2 - 2\cos(\lon-\lond_{\Nd})}\rt)\rt|_{\lon = \tilde{\lon}}
\end{bmatrix}^{T},
\end{align*}
for any $0\leq \tilde{\lon} \leq 2\pi$.  Just as $\vec{b}(\tilde{\lon})^{T}A^{-1} \uXd$ gives the value of $X(\tilde{\lon},t)$, we can use $\ubd(\tilde{\lon})$ to obtain the $n^{\text{th}}$ derivative of $X(\lon,t)$ with respect to $\lon$ as 
\begin{align*}
\lf.\dfrac{\partial^n}{\partial \lon^n} X(\lon,t)\rt|_{\lon=\tilde{\lon}}= \ubd(\tilde{\lon})^{T}A^{-1} \uXd.
\end{align*}
The evaluation of the $n^{\text{th}}$ derivative of $X(\lon,t)$ at the data site nodes $\lond_1,\ldots,\lond_{\Nd}$ can then be obtained as follows:
\begin{align}
\begin{bmatrix}
\lf.\dfrac{\partial^n}{\partial \lon^n} X(\lon,t)\rt|_{\lon=\lond_1} \\
\vdots \\
\lf.\dfrac{\partial^n}{\partial \lon^n} X(\lon,t)\rt|_{\lon=\lond_{\Nd}}
\end{bmatrix}
= 
\underbrace{
\begin{bmatrix}
\ubd(\lond_1)^{T} \\
\vdots \\
\ubd(\lond_{\Nd})^{T}
\end{bmatrix}}_{\displaystyle B_{\lond}^{n}}
A^{-1}\uXd
= 
\underbrace{B_{\lond}^{n}A^{-1}}_{\displaystyle \Dlond{n}} \uXd.
\label{eq:diff_mat_lond}
\end{align}
In a similar manner, the evaluation of the $n^{\text{th}}$ derivative of $X(\lon,t)$ at the sample site nodes $\lons_1,\ldots,\lons_{\Ns}$ can be obtained by
\begin{align}
\begin{bmatrix}
\lf.\dfrac{\partial^n}{\partial \lon^n} X(\lon,t)\rt|_{\lon=\lons_1} \\
\vdots \\
\lf.\dfrac{\partial^n}{\partial \lon^n} X(\lon,t)\rt|_{\lon=\lons_{\Ns}}
\end{bmatrix}
= 
\underbrace{
\begin{bmatrix}
\ubd(\lons_1)^{T} \\
\vdots \\
\ubd(\lons_{\Ns})^{T}
\end{bmatrix}}_{\displaystyle B_{\lons}^{n}}
A^{-1}\uXd
= 
\underbrace{B_{\lons}^{n}A^{-1}}_{\displaystyle \Dlons{n}} \uXd.
\label{eq:diff_mat_lons}
\end{align}
For given data sites $\uXd$ at any time $t$, we can interpolate these values with an RBF expansion \emph{and} evaluate the $n^{\text{th}}$ derivative of the interpolant at the data site nodes by the matrix-vector product $\Dlond{n}\uXd$ and at the sample site nodes by $\Dlons{n}\uXd$.  We refer to the $\Nd \times \Nd$ matrices $\Dlond{n}$ and the $\Ns \times \Nd$ matrices $\Dlons{n}$ as \emph{RBF differentiation matrices}.

The matrices $\Dlond{n}$ and $\Dlons{n}$ can be used to give values at respective data site or sample site nodes of the $n^{\text{th}}$ derivative of the RBF interpolant of any quantity whose values we have at the data site nodes (\emph{e.g.}, $\uYd = [\Yd_1(t)  \cdots  \Yd_{\Nd}(t)]^{T}$). These matrices can also be precomputed once at $t=0$ and used for all subsequent times.

Having defined the operators to compute derivatives of the RBF interpolant, we define the quantity
\begin{align}
\tngt := \frac{\partial}{\partial \lon} \vX(\lon,t) = \left(\frac{\partial}{\partial \lon} X(\lon,t),\frac{\partial}{\partial \lon} Y(\lon,t)\right) = (\tau_X,\tau_Y).
\label{eq:tngt_2d}
\end{align}

The unit tangent vector to $\vX(\lon,t)$ is then given by
\begin{align}
\utngt: &= \frac{\tngt}{\|\tngt\|} = (\hat{\tau}_X,\hat{\tau}_Y).  \label{eq:utngt_2d}
\end{align}
In our experiments, we assume that the Lagrangian force at a point on a platelet is the sum of a tension force, a bending-resistant force and possibly a force due to a bond between that point and a point on another platelet or the vessel wall. For the tension force,
we use the fiber model defined in~\cite{Peskin:2002}, according to which the elastic tension force density at
$\vX(\lon_i,t_k)$ is given by
\begin{align}
\frcT(\lon_i,t_k) = \left.\frac{\partial}{\partial \lon} (T\utngt)\right|_{\lon_i,t_k}, \label{eq:force_tension_2d}
\end{align}
\noindent where $T = \kt(\|\tngt\| - l_0)$ is the fiber tension and $\kt >0$ is constant. We set $l_{0,i} = \left.\|\tngt\| \right|_{\lon_i,t_0}$, where $t_0$
is the initial time of the simulation. For a bending-resistant force, we use a linear variant of the force defined in ~\cite{Griffith:2009} and define the elastic force density at $\vX(\lon_i,t_k)$ due to how much the platelet surface there is bent to be
\begin{align}
\frcB(\lon_i,t_k) = -\left.\kb\left(\frac{\partial^4 \vX}{\partial \lon^4} - \frac{\partial^4 \vX^0}{\partial \lon^4}\right)\right|_{\lon_i,t_k}. \label{eq:force_bending2_2d}
\end{align}
\noindent
Here $\vX^0 = \vX(\lon_i,t_0)$ is the initial configuration of the platelet and
$\kb > 0$ is constant. Ideally, the constants $\kt$ and $\kb$ would be chosen to reflect values determined from experiments involving real platelets. In our work, we choose $\kt$ and $\kb$ that keep isolated platelets in simple shear flows approximately rigid, and scale them as we refine the background Eulerian grid; this scaling ensures that the mechanical properties of the elastic material converge as the background grid is refined~\cite{FauciPeskin88}. Despite this approximate rigidity of isolated platelets, these platelets may deform significantly due to interactions (binding and unbinding) with other platelets, and in other models, as a consequence of platelet activation also.

We defer discussion of how we compute the forces given by Equations \eqref{eq:force_tension_2d} and \eqref{eq:force_bending2_2d} to the next
section (and the Appendix), since the implementation is different for the RBF and piecewise-linear representations of the platelet boundary. However,
the force acting on a platelet due to other platelets (and/or walls)
is common to both methods.  We use the spring force defined in~\cite{Fogelson2008}: let $p_1,p_2,...,p_{N_p}$ be the indices corresponding to the platelets in the domain. Let $p_1$ and $p_2$ be the indices of two platelets which are linked at sample sites $\vX_{p_1}(\lon^s_{i_1})$ and $\vX_{p_2}(\lon^s_{i_2})$. The force at $\vX_{p_1}(\lon^s_{i_1})$ is then given by: 
\begin{equation}
\label{eq:force_IP}
\frcC_{p_1}(\lon^s_{i_1},t_k) = \KC(||\vX_{p_2}(\lon^s_{i_2}) - \vX_{p_1}(\lon^s_{i_1})|| - \ellC)\frac{\vX_{p_2}(\lon^s_{i_2}) - \vX_{p_1}(\lon^s_{i_1})}{||\vX_{p_2}(\lon^s_{i_2}) - \vX_{p_1}(\lon^s_{i_1})||},
\end{equation}
where $\KC$ and $\ellC$ are the interplatelet cohesion spring stiffness and
the resting length, respectively; we also set $\frcC_{p_2}(\lon^s_{i_2},t_k) = - \frcC_{p_1}(\lon^s_{i_1},t_k)$. The formulation for platelet-wall links is similar. 
\section{Numerical Discretization}
\label{sec:implementation}

In this section we present the implementation details for both IB methods. For each method,
we briefly describe the spatial discretization for both the Lagrangian and Eulerian quantities. We then
describe the time-stepping scheme for each method.

\subsection{The Piecewise-Linear IB method}

Traditionally, finite-difference approximations of Equations \eqref{eq:force_tension_2d} and \eqref{eq:force_bending2_2d} are
used in conjunction with piecewise linear methods in 2D ({\em e.g.}~\cite{Griffith:2009}).
We use a second-order central difference involving sets of sample sites or IB points to discretize the derivatives
involved in the computation of both the tension and bending forces (including
tangent lengths). It is useful to think of these finite difference approximations
to the constitutive model as Hookean springs connecting pairs of IB points. Note that these differences are only second-order assuming a near-uniform sampling. This is one of the sources of error for the IB method.

For the Eulerian spatial discretization, we use a second-order centered finite-difference approximation to the Laplacian on a staggered MAC grid~\cite{HarlowWelch65}. We discretize the advection term (in conservative form $\nabla \cdot (\vu \vu^T)$ using second-order centered differences, averaging quantities to cell edges or nodes as required. For the approximate $\delta$-function, we use the ``cosine'' form
described by Peskin~\cite{Peskin:2002} which ensures that the
entire IB force is transmitted to the grid, that the force density on
the grid is a continuous function of the IB point locations, and that
the communication between grid and IB points is very localized. We note that one could also use other discrete $\delta$-functions designed to satisfy specific properties, if required~\cite{LiuMori2012}. To prevent leakage, the tension stiffness $k_t$ is set sufficiently high so that the IB point spacing on the surface is approximately $0.5 h$ (or less), where $h$ is the Eulerian grid cell width. After each update of the IB point locations, new links are formed and existing ones are broken using the model's rules for these types of events. 

We use the formally second-order Runge-Kutta time-stepping scheme outlined in \cite{DevendranPeskin}. This time-stepping scheme demonstrates second-order convergence in time for a smooth forcing function, or for an elastic material that fills the whole domain, as demonstrated in \cite{DevendranPeskin}. This scheme exhibits only first-order convergence in time in the presence of a sharp interface between the fluid and the elastic material, as is typical of IB methods. The full scheme is presented in Appendix B.

\subsection{The RBF-IB method}

In order to construct the operators utilized by
our algorithm, we must first choose an appropriate
node set. We use $\Nd$ equally-spaced values on the interval
$(0,2\pi]$ as the data site node set $\{\lond_k\}_{k=1}^{\Nd}$.
This gives a uniform sampling in the parametric space. We also use $\Ns > \Nd$ (typically, $\Ns = 4 \Nd$ or $\Ns = 8 \Nd$) equally-spaced points in the interval $(0,2\pi]$ as the set of sample site nodes
$\{\lons_j\}_{j=1}^{\Ns}$ since this results in a set of sample
sites that are well distributed over the object. As in the traditional IB method, we make sure to start simulations with a sample site spacing of less than $0.5 h$ (again enforced approximately using the tension stiffness), though the data site spacing can be much greater. In the results section, we explore the ramifications of this choice.

We have formulated our operators to ensure that operations like evaluation
of the interpolant and computing derivatives (and therefore the
constitutive model) do not require solving a linear system for any time step of the platelet simulation except the initial step.  This is possible because, though the data sites and sample sites move over the course of the simulation, their values in parameter space do not change. For the RBF model of the platelets, the evaluation matrix $\Es$ in Equation \eqref{eq:rbf_eval} and differentiation matrices $\Dlond{n}$ and $\Dlons{n}$ in Equations \eqref{eq:diff_mat_lond} and \eqref{eq:diff_mat_lons}, respectively can be computed using the FFT as discussed in our previous work~\cite{Shankar2012}.  This is possible since the data site nodes $\{\lond_k\}_{k=1}^{\Nd}$ are equally-spaced, which results in the $A$ matrix defined Equation \eqref{eq:rbf_linsys} having a circulant matrix structure. The costs and accuracy of the RBF models are elaborated upon in the discussion of the results. The algorithm to compute forces on platelets using these operators is presented in the Appendix A.

The RBF-IB method uses the same time-stepping scheme and Eulerian discretization as the piecewise linear IB method, with one important difference. When computing the forces at time level $n+1/2$, we advance the \emph{data sites} to time level $n+1/2$, generate a set of sample sites at that time level, and compute forces at the sample sites. Similarly, we use the mid-step approximation to the velocity field to advance the \emph{data sites} to time level $n+1$. We thus generate only a single set of sample sites every time-step, since the sample sites are only needed when the data sites are advanced to time level $n+1/2$. It is clear that if the number of data sites is fewer than the number of sample sites, this results in improved computational efficiency over the piecewise linear IB method. However, it is important to explore the effect of our changes on the convergence of the algorithm. We explore these questions in the results section. For a more complete description of the RBF-IB time-stepping scheme, see Appendix C.

\section{Results}
\label{sec:results}
In this section, we first compare the convergence of the RBF-IB method on a canonical test problem. We also use this test problem to explore
the relationship between the number of data sites ($\Nd$) and the Eulerian grid spacing ($h$). We then
compare the area loss in an elastic object simulated by each method on the same problem, and discuss the time-step sizes allowed by both methods. We follow
with a discussion of the change in energy over time in the RBF-IB method. We then provide timings for platelet simulations and discuss both
foreseen and unforeseen results of using the RBF model within the IB method. Finally, we present the results of platelet aggregation simulations conducted using the RBF-IB method.

\textbf{Description of our standard fluid-structure interaction problem:}

\begin{figure}[htbp]
\centering
\includegraphics[height=3.0in,width=3.0in]{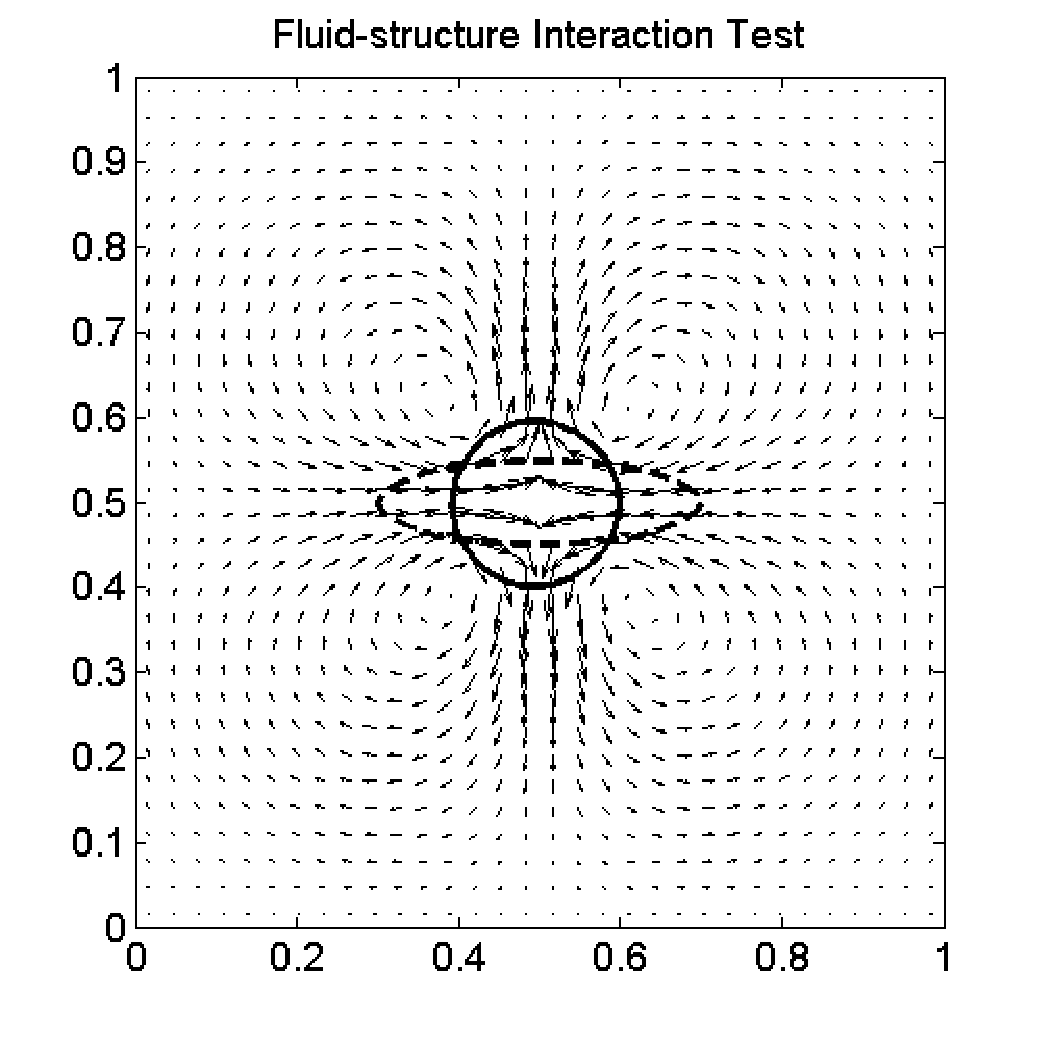}
\caption{A visualization of the fluid-structure interaction test. The dashed lines show the initial ellipse, while the filled line indicates the near-circular object at the final time $t = 2.0$. The small
arrows indicate the velocity field at the final time. The maximum velocity is very close to zero at this time as the object is almost at rest.}
\label{fig:elastest}
\end{figure}

We describe a standard fluid-structure interaction problem on which we test both versions of the IB method. This problem is commonly used in the IB literature (\emph{e.g.}, \cite{Newren2007}). The problem involves placing an elliptical object with its center of mass at the center of the $[0,1]^2$ physical domain. The elliptical object has a circle of the same area as its rest configuration, and attempts to attain the rest configuration subject to a combination of tension and bending forces. The physical domain is filled with a fluid that is initially at rest, with periodic boundary conditions in the $x$-direction and no-slip Dirichlet boundary conditions in the $y$-direction. We set the radius of the target circle to be $r = 0.1$ units, with the ellipse initially having a major axis of $a = 2r$ and a minor axis of $b = 0.5r$. This test is visualized in Figure \ref{fig:elastest}.

\subsection{Convergence studies}

In previous work~\cite{Shankar2012}, we compared the accuracy and convergence of both RBF and traditional IB geometric modeling strategies for static platelet-like shapes. We now
compare the accuracy and convergence of the full RBF-IB and PL-IB methods, both for the velocity field and for the immersed elastic structure.

For the fluid, on each grid with cell width $h$, we define the quantity ${\vu}^{c,h}$, the coarsened discrete velocity field from the $256\times256$ grid. This coarsened velocity is obtained using a spline interpolation of each component of the velocity on the finest grid, and evaluating it at the edges on the coarser grid. For each grid point $i$ on a grid with cell width $h$, we compute the quantity $f^h_i = ||{\vu}^h_i - {\vu}^{c,h}_i||$. We define the $l_2$ error in the velocity field as $e_{2}(h) = \sqrt{\sum_i {(f^h_i h)}^2}$, and the $l_{\infty}$ error in the velocity field to be $e_{\infty}(h) = \max_i f^h_i$. The convergence rate for errors $e(2h)$ and $e(h)$ is measured as $p_u = \log_2\lf(\frac{e(2h)}{e(h)}\rt)$.

We note that convergence rates of numerical solutions can be computed in an alternate manner to the one outlined above. Solutions on three grids $\vu^{4h}$, $\vu^{2h}$ and $\vu^{h}$ are chosen; the errors $e(2h)$ and $e(h)$ are computed by comparing the $4h$ grid to the $2h$ grid, and the $2h$ grid to the $h$ grid respectively. The rate of convergence is then once again estimated as $p_u = \log_2\lf(\frac{e(2h)}{e(h)}\rt)$. While the results shown in this section do not employ this method, we have verified that the orders of convergence computed using this method match those presented in this section.

For the Lagrangian markers (sample sites or IB points), we adopt the following procedure:
\begin{enumerate}
\item Given the number of sample sites $\Ns$ and the radius of the target circle $r$, we define $\theta = \frac{2\pi}{\Ns}$, the angle subtended at the center of the circle if the points were evenly-spaced.
\item We define the quantity $C = 2r\sin(0.5\theta)$, the chord length between any two points in a set of evenly-spaced points on an ideal circle. We also define $C_{exact}$ to be the ideal chord length for $\Ns = 400$.
\item We compute the actual distances $d_i$ between the sample sites (or IB points) for a simulation computed on the $256\times256$ grid with $\Ns = 400$. We then compute the quantities $s_{\infty} = \max_i |d_i - C_{exact}|$ and $s_{2} = (1/{\Ns})\sqrt{\sum_i |d_i - C_{exact}|^2}$.
\item We compute $s^{\Ns}_{\infty}$ and $s^{\Ns}_{2}$ for $\Ns = 50,100,200$. We then define the $l_2$ error to be $e_{2}(\Ns) = |s^{\Ns}_2 - s_{2}|$ and the $l_{\infty}$ error to be $e_{\infty}(\Ns) = |s^{\Ns}_{\infty} - s_{\infty}|$. 
\end{enumerate}
We define the convergence rate for errors $e(\Ns)$ and $e(2\Ns)$ to be $p_{X} = \log_2\lf(\frac{e(\Ns)}{e(2\Ns)}\rt)$. 

By defining the Lagrangian errors in this fashion, we circumvent the fact that the `correct' sample site spacing is unknown; our definitions above measure `errors' against equally-spaced points on the circle. However, as can be seen in the following text, these errors converge at a first-order rate (asymptotically). For full transparency, we also record the deviations of the sample sites from an equispaced set of IB points on the finest Lagrangian `grid'; these are given by the quantities $s_{\infty}$ and $s_2$.

\subsubsection{Convergence of the RBF-IB method for $\Nd = 25$:}

\begin{table}[ht]
\begin{tabular}{|c|c|c|c|c|c|c|}
\hline
Grid Size & {$\Ns$} & $\Delta t$ & $L_2$ error & Order of convergence & $L_{\infty}$ error & Order of convergence \\ 
\hline
{$32\times32$} & 50 & {$2\times10^{-4}$} & 1.7666e-04 & & 8.8682e-04  &  \\\hline
{$64\times64$} & 100 & {$1\times10^{-4}$} & 1.5097e-04 & 0.23 & 5.4255e-04 & 0.71\\\hline 
{$128\times128$} & 200 & {$5\times10^{-5}$} & 8.4247e-05 & 0.84 & 3.0738e-04 & 0.82 \\\hline
\end{tabular}	
\caption{Results of a refinement study with the RBF-IB method with $\Nd = 25$ data sites. We show the convergence of the velocity field, with errors measured against
the velocity field of a simulation on a $256\times256$ grid with $\Ns = 400$ sample sites, $\Nd = 100$ data sites and $\Delta t = 2.5\times10^{-5}$.}	
\label{tab:RBFIBF1}
\end{table}

\begin{table}[ht]
\begin{tabular}{|c|c|c|c|c|c|c|}
\hline
{$\Ns$} & Grid Size & $\Delta t$ & $L_2$ error & Order of convergence & $L_{\infty}$ error & Order of convergence \\
\hline
{$50$}  & {$32\times32$} & {$2\times10^{-4}$} & 3.1188e-06 &  & 2.3238e-05  &  \\\hline 
{$100$} &  {$64\times64$} & {$1\times10^{-4}$} & 3.5898e-07 & 3.12   & 3.0048e-06  & 2.95 \\\hline 
{$200$} &  {$128\times128$} & {$5\times10^{-5}$} & 1.5310e-07 & 1.23   & 1.7905e-06 & 0.75 \\\hline
\end{tabular}	
\caption{Results of a refinement study with the RBF-IB method with $\Nd = 25$ data sites. We show the convergence in the sample site positions, with errors measured against
the sample site positions of a simulation on a $256\times256$ grid with $\Ns = 400$ sample sites, $\Nd = 100$ data sites and $\Delta t = 2.5\times10^{-5}$.}	
\label{tab:RBFIBP1}
\end{table}

Here, we test the convergence of the RBF-IB method on the fluid-structure interaction problem described above. We compare the velocity field and sample site positions to those computed for the same test problem on a $256\times256$ grid with $\Ns = 400$ sample sites and $\Nd = 100$ data sites. Tables \ref{tab:RBFIBF1} and \ref{tab:RBFIBP1} show the results obtained with $\Nd = 25$ data sites for the fluid and the object respectively. 

Table \ref{tab:RBFIBF1} shows that the errors produced in the RBF-IB method. We see that the RBF-IB method for $\Nd = 25$ shows first-order convergence asymptotically. Table \ref{tab:RBFIBP1} shows results for the errors on the object. The errors on the coarsest grid are high, leading to a higher-than-expected convergence rate in both norms when we measure the errors on a $64\times64$ grid. Once again, the convergence rate remains close to first-order, as expected. For completion, we note that $s_2 = 4.3694e-08$ and $s_{\infty} = 1.1113e-06$ for the structure on the $256\times256$ grid for $\Nd = 25$ data sites.

\subsubsection{Convergence of the RBF-IB method for $\Nd = 50$:}

\begin{table}[ht]
\begin{tabular}{|c|c|c|c|c|c|c|}
\hline
Grid Size & {$\Ns$} & $\Delta t$ & $L_2$ error & Order of convergence & $L_{\infty}$ error & Order of convergence \\ 
\hline
{$32\times32$} & 50 & {$2\times10^{-4}$} & 5.5617e-03 & & 3.7404e-02  &  \\\hline
{$64\times64$} & 100 & {$1\times10^{-4}$} & 2.2436e-04 & 4.63 & 1.2403e-03 & 4.91\\\hline 
{$128\times128$} & 200 & {$5\times10^{-5}$} & 7.8934e-05 & 1.51 & 2.8859e-04 & 2.10 \\\hline
\end{tabular}	
\caption{Results of a refinement study with the RBF-IB method with $\Nd = 50$ data sites. We show the convergence of the velocity field, with errors measured against
the velocity field of a simulation on a $256\times256$ grid with $\Ns = 400$ sample sites, $\Nd = 100$ data sites and $\Delta t = 2.5\times10^{-5}$.}	
\label{tab:RBFIBF2}
\end{table}

\begin{table}[ht]
\begin{tabular}{|c|c|c|c|c|c|c|}
\hline
{$\Ns$} & Grid Size & $\Delta t$ & $L_2$ error & Order of convergence & $L_{\infty}$ error & Order of convergence \\
\hline
{$50$}  & {$32\times32$} & {$2\times10^{-4}$} & 3.6212e-05 &   & 5.2794e-04  &  \\\hline 
{$100$} &  {$64\times64$} & {$1\times10^{-4}$} & 3.8824e-07 & 6.54  & 7.8472e-06  & 6.07 \\\hline 
{$200$} &  {$128\times128$} & {$5\times10^{-5}$} & 1.7236e-07 & 1.17   & 2.0636e-06 & 1.93 \\\hline
\end{tabular}	
\caption{Results of a refinement study with the RBF-IB method with $\Nd = 50$ data sites. We show the convergence in the sample site positions, with errors measured against
the sample site positions of a simulation on a $256\times256$ grid with $\Ns = 400$ sample sites, $\Nd = 100$ data sites and $\Delta t = 2.5\times10^{-5}$.}	
\label{tab:RBFIBP2}
\end{table}

We repeat the above test problem with $\Nd = 50$ data sites. As before, we compare the velocity field and sample site positions to those computed for the same test problem on a $256\times256$ grid with $\Ns = 400$ sample sites and $\Nd = 100$ data sites. Tables \ref{tab:RBFIBF2} and \ref{tab:RBFIBP2} show the results for the fluid and the object respectively. 

Examining Table \ref{tab:RBFIBF2}, we see that when moving to the finest grid, we have once again recovered first-order convergence. The errors on the $128\times128$ grid for $\Nd = 50$ with the RBF-IB method are close to those on the same grid with $\Nd = 25$. Table \ref{tab:RBFIBP2} shows errors similar to those seen in Table \ref{tab:RBFIBP1}, albeit with less erratic convergence. Indeed, we recover first-order convergence in the $l_2$ norm and close to second-order convergence in the $l_{\infty}$ norm.

We note that using $\Nd = 50$ data sites does not result in significantly better convergence on the structure than $\Nd = 25$. There are two possible explanations. The first is that the function representing the shape of the object is of limited smoothness (as seen in our previous work \cite{Shankar2012}), with higher values of $\Nd$ causing the interpolation error to saturate or even increase. The alternate (and more likely) explanation is that, since our RBFs are parametrized on the circle, $\Nd = 25$ would already have a very high accuracy when the object becomes a circle, considering the spectral accuracy of RBF interpolation on the circle; in such a scenario, using $\Nd = 50$ data sites would only serve to increase the rounding errors in the representation of the structure. The values of $s_2$ and $s_{\infty}$for the structure on the $256\times256$ grid for $\Nd = 50$ data sites are the same as those for $\Nd = 25$ data sites.

\subsubsection{Convergence of the RBF-IB method for $\Nd = 0.25\Ns$:}

\begin{table}[ht]
\begin{tabular}{|c|c|c|c|c|c|c|}
\hline
Grid Size & {$\Ns$} & $\Delta t$ & $L_2$ error & Order of convergence & $L_{\infty}$ error & Order of convergence \\ 
\hline
{$32\times32$} & 50 & {$2\times10^{-4}$} & 4.7909e-04 & & 2.4700e-03  &  \\\hline
{$64\times64$} & 100 & {$1\times10^{-4}$} & 1.5097e-04 & 1.67 & 5.4255e-04 & 2.19\\\hline 
{$128\times128$} & 200 & {$5\times10^{-5}$} & 9.0802e-05 & 0.73 & 3.3020e-04 & 0.72 \\\hline
\end{tabular}	
\caption{Results of a refinement study with the RBF-IB method with $\Nd = 0.25\Ns$ data sites. We show the convergence of the velocity field, with errors measured against
the velocity field of a simulation on a $256\times256$ grid with $\Ns = 400$ sample sites, $\Nd = 100$ data sites and $\Delta t = 2.5\times10^{-5}$.}	
\label{tab:RBFIBF3}
\end{table}

\begin{table}[ht]
\begin{tabular}{|c|c|c|c|c|c|c|}
\hline
{$\Ns$} & Grid Size & $\Delta t$ & $L_2$ error & Order of convergence & $L_{\infty}$ error & Order of convergence \\
\hline
{$50$}  & {$32\times32$} & {$2\times10^{-4}$} & 9.7439e-06 &   & 7.0459e-05  &  \\\hline 
{$100$} &  {$64\times64$} & {$1\times10^{-4}$} & 3.5898e-07 & 4.76   & 3.0048e-06  & 4.55 \\\hline 
{$200$} &  {$128\times128$} & {$5\times10^{-5}$} & 1.2694e-07 & 1.50   & 1.4362e-06 & 1.07 \\\hline
\end{tabular}	
\caption{Results of a refinement study with the RBF-IB method. We show the convergence in the sample site positions, with errors measured against
the sample site positions of a simulation on a $256\times256$ grid with $\Ns = 400$ sample sites, $\Nd = 100$ data sites and $\Delta t = 2.5\times10^{-5}$.}	
\label{tab:RBFIBP3}
\end{table}

In the traditional IB method, the number of IB points depends on the grid spacing $h$. Typically, the number of IB points is chosen
so that the distance between any two sample sites is always less that $0.5h$. In all the tests above, we have maintained that relationship
for the sample sites in the RBF-IB method. In the RBF-IB method, we always use fewer data sites than sample sites, \emph{i.e.}, $\Nd < \Ns$, with the choice of $\Nd$ being justified by the results in previous work~\cite{Shankar2012}. Furthermore, in the tests above, we fix $\Nd$ even as we refine the fluid grid. For $\Nd = 25$, this means that as we refine $\Ns$ the distance between data sites increases from $0.8h$ to $3.2h$ (at the start of the simulation).

In order to gain intuition on the relationship between $\Nd$ and $h$, we now perform a convergence study (using the same test problem given above) with increasing values of $\Nd$ as $h$ is reduced. To accomplish this, we use values of $\Nd = 12, 25, 50,100$ for $\Ns = 50,100,200,400$, \emph{i.e.}, we enforce $\Nd = 0.25\Ns$. We use the solution computed with $\Nd = 100$ and $\Ns = 400$ on a $256\times256$ grid as our gold standard, just as we have in all the other tests. Table \ref{tab:RBFIBF3} shows the results for the fluid. Clearly, the errors are higher and the convergence more erratic than for the fixed $\Nd = 50$ tests previously presented, but varying $\Nd$ certainly seems to give better convergence than fixing it to $\Nd = 25$. However, the convergence in the structure is comparable, with lower errors being achieved than both $\Nd = 25$ and $\Nd = 50$. This can be seen in Table \ref{tab:RBFIBP3}. Unfortunately, the advantages of varying $\Nd$ with $\Ns$ are not clear. Using $\Nd = 50$ yields the lowest errors in the fluid on the finest grid level, and reasonably low errors on the structure for all grid levels. Given the similarity of the errors achieved with $\Nd = 50$ to those achieved with increasing $\Nd$, we choose the simpler strategy of using a fixed value of $\Nd = 50$ for our tests, though we present timings with $\Nd = 25$ as well. 

\subsubsection{Effect of the shape parameter $\ep$:}

In previous work \cite{Shankar2012}, we found that the RBF shape parameter $\ep > 0$ had to be selected carefully to achieve spectral accuracy in the representation of the elastic structure. In that work, we found that small values of $\ep$ were ideal for interpolating smooth target shapes and larger ones for rougher target shapes. In our tests, we found that the errors depended on $\ep$ even in the case of fluid-structure interaction, with smaller values of $\ep$ giving the lowest values of $s^e_2$ and $s^e_{\infty}$ on the $256\times256$ grid. However, as we mentioned in our previous work, lower values of $\ep$ can make the RBF interpolation matrix more ill-conditioned. While methods (such as RBF-QR and RBF-RA) have been developed to overcome this poor conditioning~\cite{FornbergPiret:2007}, they are much more expensive than forming and inverting the standard RBF interpolation matrix. We thus choose a small value of $\ep = 1.2$ for all our tests. When using $\Nd = 100$, we use $\ep = 2.0$ (which was verified on a static test case to be accurate to 12 digits). These are the smallest values we were able to pick without the interpolation matrix becoming ill-conditioned, a strategy consistent with the one used in our previous work \cite{Shankar2012}. 

\subsection{Area loss and time-step size in the RBF-IB method:}

In this section, we study the area loss in the RBF-IB method in a refinement study. We then explore the maximum stable time-step size afforded by each IB method.

\begin{table}[ht]
\begin{tabular}{|c|c|c|c|c|c|}
\hline
{$\Ns$} & Grid Size & $\Delta t$ & $\%$ area loss ($\Nd = 25$) & $\%$ area loss ($\Nd = 50$) & $\%$ area loss ($\Nd = \Ns/4$) \\
\hline
{$50$}  &  {$32\times32$}   & {$2\times10^{-4}$}  & 0.0680  & 0.3081  &  0.0450 \\\hline 
{$100$} &  {$64\times64$}   & {$1\times10^{-4}$} & 0.0047  & 0.0049  &  0.0047     \\\hline 
{$200$} &  {$128\times128$} & {$5\times10^{-5}$} & 0.0023  & 0.0025 &  0.0025     \\\hline
{$400$} &  {$256\times256$} & {$2.5\times10^{-5}$}& 0.0015  & 0.0015 &  0.0015     \\\hline
\end{tabular}	
\caption{Percentage area loss in the RBF-IB method as a function of grid size, the number of sample sites $\Ns$ and the time step $\Delta t$. The PL-IB method
gives area losses similar to the $\Nd = 50$ case, except on the coarsest grid, where the percentage area loss is three times that of the RBF-IB method.}	
\label{tab:RBFIBAL}
\end{table}

The PL-IB method generally attempts to maintain an IB point separation distance of $0.5h$ in order to reduce area loss over the coarse of the simulation. In the RBF-IB method, while the sample
site spacing is initially set at $0.5h$, we initialize the structure with a much coarser data site discretization, with the data site separation being almost $3.2 h$ in some cases. In addition, we use
the same strategy for interpolating velocities that we do in the PL-IB method, \emph{i.e.}, we interpolate velocities to data sites from a $4\times4$ patch of fluid around each data site. While
this can result in significant computational savings, it is important to explore the area loss in our discretization. We turn once again to our standard fluid-structure interaction problem. We run that simulation on successively finer grids until time $t=2$. For both the RBF-IB method and the PL-IB method, we measure the initial area of the object for the same initial configuration of points. We then measure the area at time $t=2$ and compute the percentage change in area.

In order to get an accurate estimate of the area in both methods, we fit an RBF interpolant to each object's Lagrangian markers (data sites for the RBF-IB method and all the IB points for the PL-IB method). We then sample that interpolant at a fixed number of points ($400$ points), and use the trapezoidal rule to compute the area. As was mentioned earlier, we ensure that the initial ellipse has the same area as the target circle by picking its radii to be $a = 2r$ and $b = 0.5r$, where $r = 0.1$ is the radius of the target circle. The exact area is then $\frac{\pi}{100}$. Our approach of sampling
each object and computing the area with the trapezoidal rule gives an area estimate that agrees with this value up to 7 digits at $t=0$. We record the results of our refinement study in Table \ref{tab:RBFIBAL}. 

From the table, it is clear that the area loss for fixed $\Nd = 25, 50$ and $\Nd = 0.25\Ns$ are all close to each other. On the coarsest grid, it appears that smaller values of $\Nd$ result in lower area loss. The area losses for $\Nd = 50$ match with those given by the PL-IB method (results not shown), except in the case of the coarsest grid, where the PL-IB method gives almost a $1\%$ area loss. The convergence is initially second-order but quickly saturates. This saturation is likely due to two sources of error: the first is the interpolation of velocities to the Lagrangian markers, which does not preserve the divergence-free nature of the fluid velocity; the second is the fact that the time-integration itself is not specifically designed to preserve area. Nevertheless, it is clear from this study that the RBF-IB method produces similar area losses to the PL-IB method despite using a smaller number of Lagrangian markers to move the structure through the fluid.

Another measure of interest is the maximum stable time-step size afforded by each method. We measure this by increasing the time-step size in small increments and observing the forces produced on the structure in the fluid-structure interaction test. Using a time-step that is too large can result in the forces blowing up and the simulation halting. We immediately note that the PL-IB method allows a maximum time-step size of $\Delta t = 2\times10^{-4}$ on the $32\times32$ grid when $\Ns = 50$ IB points are used, and a maximum time-step size of $\Delta t = 10^{-4}$ on the $64\times 64$ grid when $\Ns = 100$ IB points are used. We use these values of $\Delta t$ as the starting point when testing for the time-step sizes allowed by the RBF-IB method, and increase the value of $\Delta t$ in increments of $10^{-4}$. We found that on the $32\times32$ grid, the RBF-IB method allows us to take time-steps that are $3\times$ larger than the time-steps allowed by the traditional IB method; on the $64\times64$ grid, the RBF-IB method can use time-steps that are $1.5\times$ larger than the time-steps allowed by the traditional IB method. This pattern holds both when $\Nd = 25$ and $\Nd = 50$ data sites are used. 

In simulations involving platelet-like shapes (ellipses that attempt to maintain their elliptical configuration), we found that the RBF-IB method allows time-step sizes that are $6\times$ larger than those allowed by the PL-IB method on a $32\times32$ grid, and $3\times$ larger than those allowed by the PL-IB method on a $64\times64$ grid. This is likely due to the fact that platelet simulations involve smaller deformations than those seen in the standard fluid-structure interaction test.

\subsection{Energy Estimates}

\begin{figure}[htbp]
\centering
\includegraphics[height=3.0in,width=3.0in]{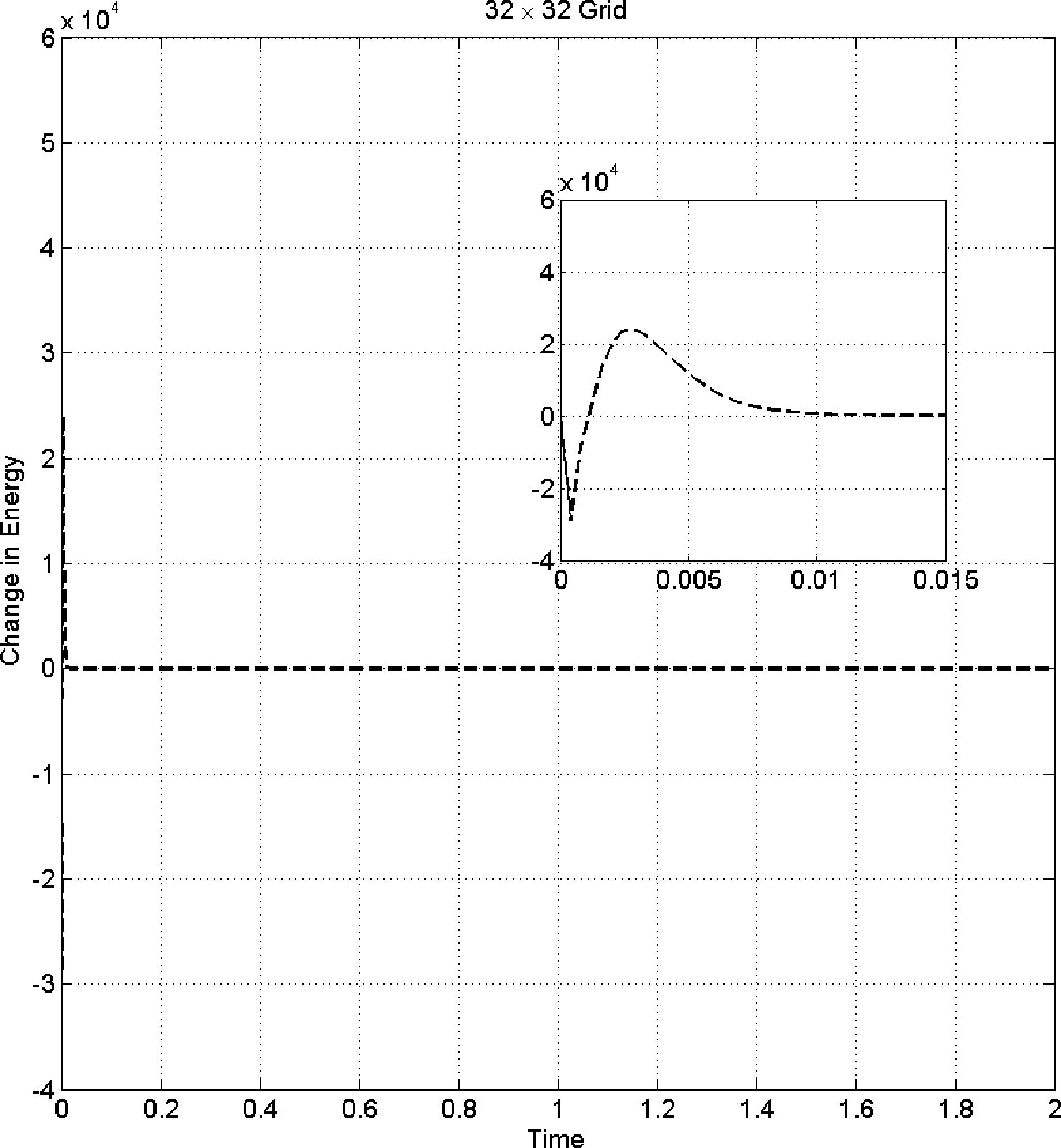}
\includegraphics[height=3.0in,width=3.0in]{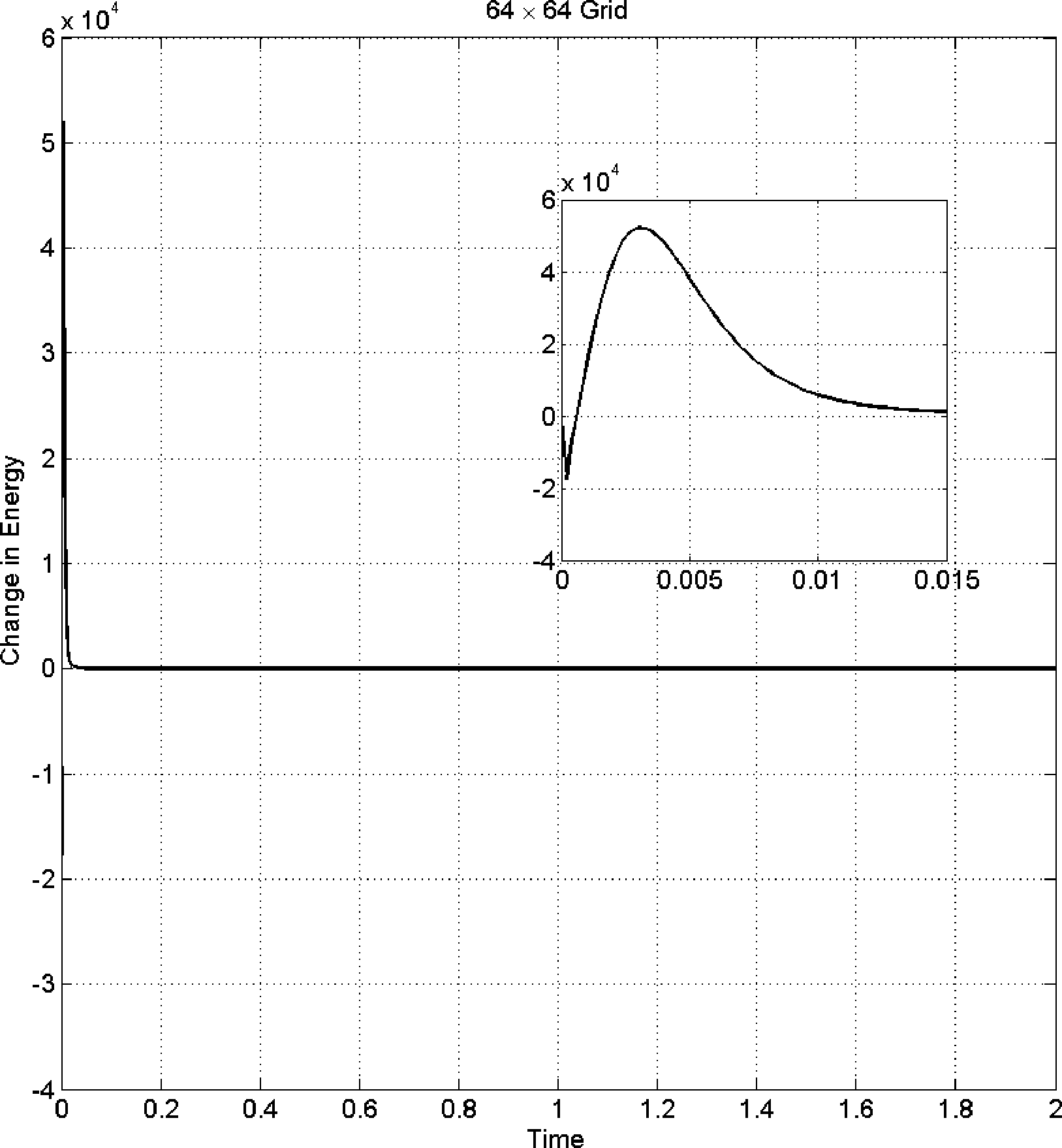}
\caption{Change in energy \emph{per time-step} as a function of time in the RBF-IB method. The figure on the top left shows the change in energy over a time-step as a function of time for the standard fluid-structure interaction test on a $32\times32$ grid. The figure on the top right shows the same quantity on a $64\times64$ grid. We use $\Nd = 50$ data sites for both grid sizes. The inset plots show the initial
spikes corresponding to the change from an ellipse to a circle which are difficult to see in the main plots.}
\label{fig:energy}
\end{figure}

In this section, we compute energy estimates for the RBF-IB method in the context of our standard
fluid-structure interaction problem. We run our simulation out to time $t = 2.0$ on two grid sizes, $32\times32$
and $64\times64$, with time-step sizes $\Delta t = 2\times10^{-4}$ and $\Delta t = 10^{-4}$ respectively. We use $\Nd = 50$ data sites.

In this test, one expects the changes in energy
to be mainly due to the deformation of the stiff elastic object. Eventually, the energy of the system
must damp out as the elastic object reaches its target configuration. We compute the change in energy to demonstrate that
the energy is bounded within the RBF-IB simulation. The energy change in a time-step is computed as the sum of the difference in
kinetic energy of the fluid over the time-step and the change in potential energy of the elastic object. This can be written as

\begin{equation}
\label{eq:rbf_energy}
\Delta E = \sum\limits_{fluid} \rho \vu^{n+1}\cdot\vu^{n+1} - \sum\limits_{fluid} \rho \vu^{n}\cdot\vu^{n} + \Delta t\sum\limits_{\vX} \vF^{n+1/2} \cdot \frac{\partial \vX}{\partial t}^{n+1/2}
\end{equation} 

Here, the Lagrangian force $\vF$ is computed at time level $n+1/2$ at the sample sites. The $\frac{\partial \vX}{\partial t}$ term is computed by
applying the evaluation matrix $\Es$ to the velocities obtained at the data sites. This gives us sample site velocities, allowing
us to compute dot products with the $\vF$ terms.

The results of this test are shown in Figure \ref{fig:energy}. Both plots show the change in energy of the system for the fluid-structure interaction problem on a $32\times32$ grid (left) and a $64\times64$ grid (right). Here, the fluid starts off stationary, so the initial
kinetic energy is zero. However, the elliptical elastic object starts off under tension, since its
target configuration is a circle. This means that the initial elastic potential energy of the system is
high (though negative by convention). As the elastic object attempts to minimize its elastic
potential energy, its deformation drives a change in the kinetic energy of the fluid, causing the kinetic
energy of the fluid to increase from its initial value of zero to some maximum. However, the elastic object
soon attains something close to its reference configuration, causing the kinetic energy of the fluid to drop sharply. The spikes in both the left and right sides of Figure \ref{fig:energy} correspond to that rise and fall in kinetic energy and the trending of the potential energy to zero on both grids, and can be seen more clearly in the inset plots. The
viscosity of the fluid causes the kinetic energy to eventually damp out almost completely, with minor perturbations due to possible deformations of the elastic object. The energy of the system
continues to decrease as the object becomes more and more circular. In fact, our estimates show that the change in energy
is negative, indicating that our method is dissipative. The results are similar for $\Nd = 25$ (not shown), though using more data sites appears to make our method less dissipative on this particular test problem.

\subsection{Timings for Platelet simulations}

We now present timings of simulations involving platelet-like shapes. The setup here is different from the standard fluid-structure interaction test. We place
ellipses ($r = 0.05$, $a=2r$, $b=0.5r$) at the left end of a $[0,2]\times[0,1]$ domain that resembles a channel. These ellipses represent platelets, and they attempt to maintain their
elliptical shapes, \emph{i.e.}, their configuration at time $t=0$ is their preferred configuration. We apply a background force that would result in parabolic velocity field in the absence of the platelets, with a density $\rho = 1.0$ and a non-dimensionalized viscosity of $\mu = 8.0$. The field has a maximum velocity value of $u_{max} = 5.0$, with no-slip boundary conditions on the top and bottom of
the domain and periodic boundary conditions at the left and right ends. A platelet is removed from the domain if its center of mass crosses the location $x = 1.9$.

\begin{figure}[htbp]
\centering
\includegraphics[height=3.0in,width=3.0in]{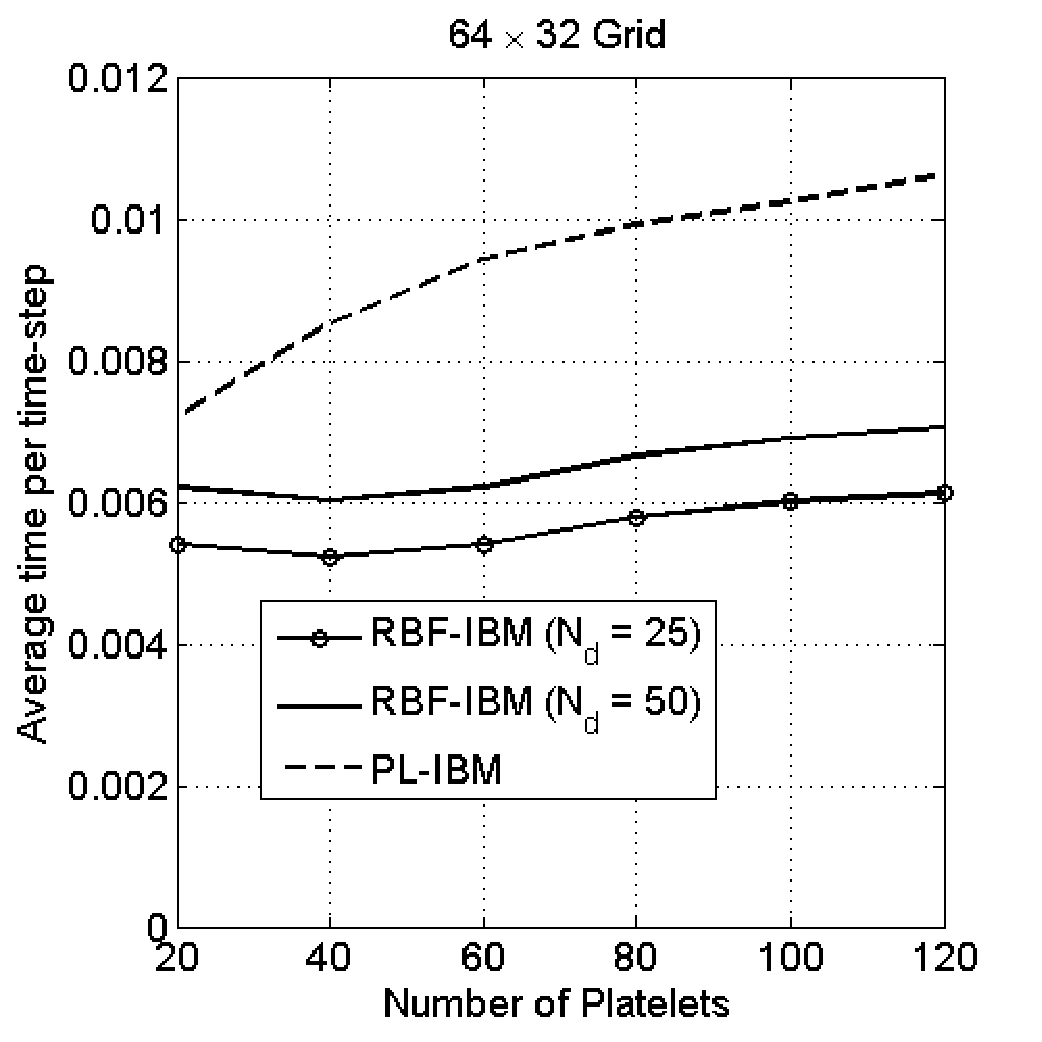}
\includegraphics[height=3.0in,width=3.0in]{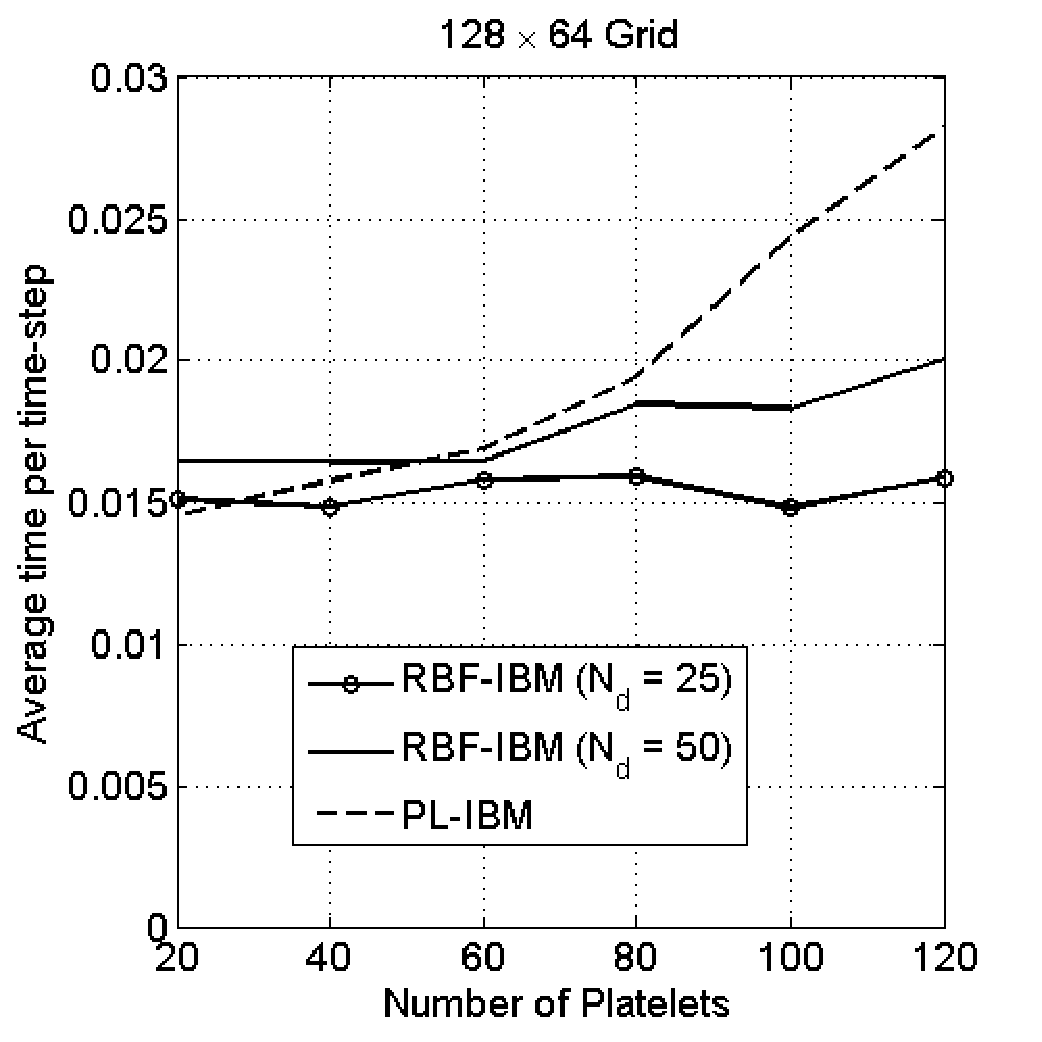}
\includegraphics[height=3.0in,width=3.0in]{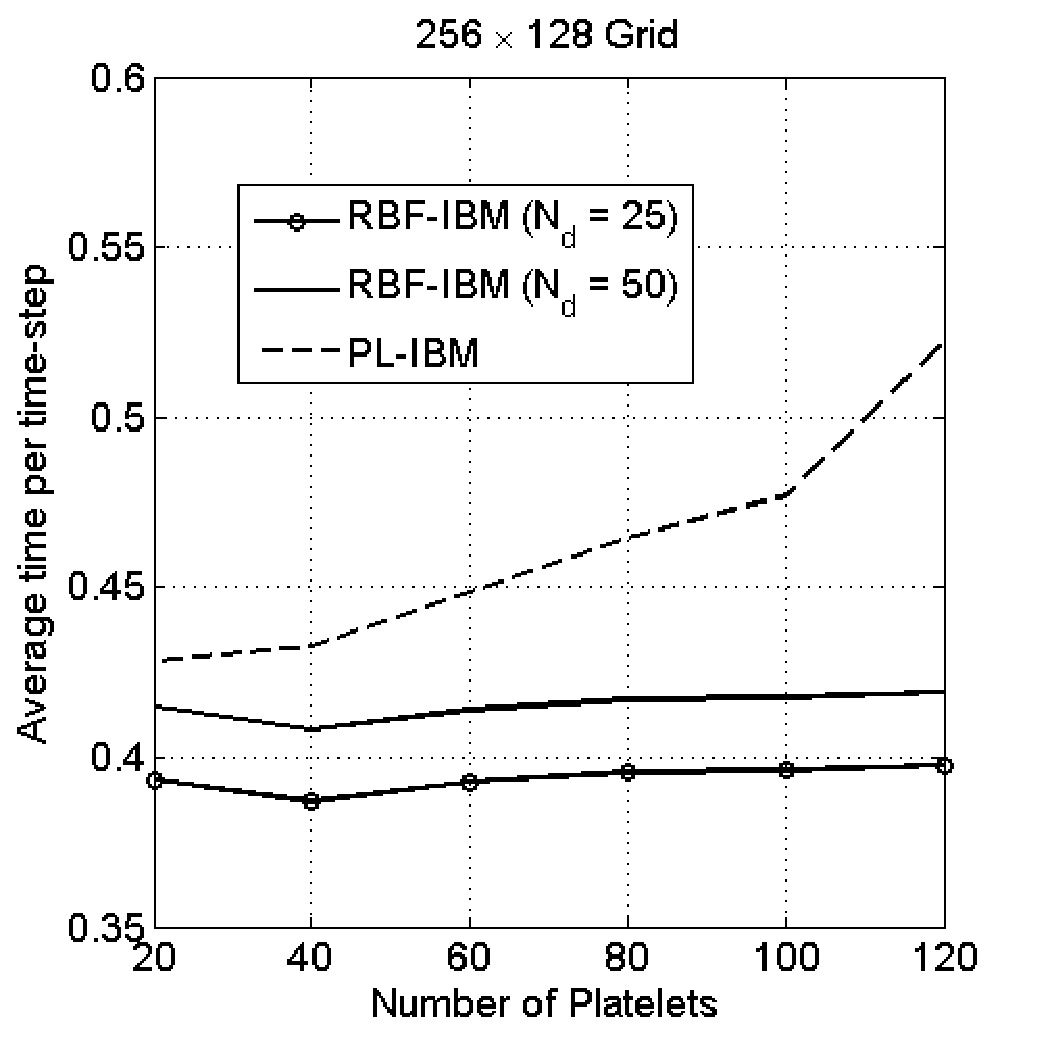}
\caption{Average time per time-step for $10^5$ time-steps of each simulation method
				as a function of the number of platelets. In the first row, the figure on the left
				shows timings on a $64 \times 32$ grid and the figure on the right on a $128 \times 64$ grid. The figure below shows
				timings on a $256 \times 128$ grid. The time-step was set to $\Delta t = 10^{-4}$
				for the figures on the top row, and was set to $\Delta t = 10^{-5}$ for
				the figure on the bottom.}
\label{fig:t1}
\end{figure}
%\FloatBarrier

Figure \ref{fig:t1} shows timings for three grid sizes for each method as a function of the number of platelets ($N_{\rm p}$)
being simulated. The number of sample sites was fixed at $\Ns = 100$ for both methods and the number of data sites
for the RBF-IB method was set to $\Nd = 25$ for one set of tests and then to $\Nd = 50$ for the next set. We plot the average time per time-step as a function of the number of platelets; this was computed by running simulations on each grid for $10^5$ time-steps, and dividing the total wall-clock time by the number of time-steps. We average this over three runs of each simulation.

While the cost of platelet operations always increases as we increase the number of platelets, the increase in cost is slower (and the absolute cost) for the RBF-IB method
due to the RBF representation. For example, for $N_{\rm p} = 60$, the PL-IB method directly spreads forces from, interpolates to, and moves a total of $6000$ IB points (twice per time-step due to the RK2 scheme) while the RBF-IB method with $\Nd = 25$ data sites computes forces at $6000$ points, and interpolates velocities to (and moves) only $1620$ points twice per time-step. If the number of platelets is doubled to $N_{\rm p} = 120$, the PL-IB method now computes forces at, spreads from, interpolates to, and moves $12000$ points twice per time-step,
whereas the RBF-IB method computes forces at $12000$ points, but interpolates velocities to and moves only $3240$ points twice per time-step. The cost of the RBF-IB method shows better than linear scaling
with respect to the number of platelets on all the tested grid sizes for these reasons. Furthermore, it is clear that there is not much of a difference in computational cost between using $\Nd = 25$ data sites and $\Nd = 50$ data sites. 

\subsubsection{Effect of the RBF representation on the fluid solver}

In previous work \cite{Shankar2012}, we showed that using an RBF interpolant for geometric modeling is more computationally efficient (for a given accuracy) than using piecewise quadratics and finite differences. However, that benefit alone does not explain the computational efficiency of the RBF-IB method over the PL-IB method that we see in Figure \ref{fig:t1}.

To fully understand the speedup seen with the RBF-IB method, it is important to understand how the costs are distributed between the different operations (platelet operations and fluid solves) in both IB methods. We show the results for $N_{\rm p} = 60$ platelets in Table \ref{tab:Profile}. Clearly, as $h$ is reduced, both IB codes spend more time in the fluid solver than on platelet operations. However, the RBF-IB method clearly spends less time in the fluid solver than the PL-IB method does as we refine the background Eulerian grid.

Indeed, this unexpected result is what gives the RBF-IB method an edge even when the cost of fluid solves dominates the cost of platelet operations. We hypothesize that this may be caused by the RBF representation producing smoother Lagrangian forces than the finite difference model used in the PL-IB method. Our experiments show that the RBF-IB code needs fewer iterations in the linear solver used in the pressure projection-- anywhere from $10-30 \%$ fewer than the identical fluid solver used in the PL-IB method, depending on the time-step size and the grid resolution, with larger savings on finer grids and smaller time-step sizes.

\begin{table}[ht]
\begin{tabular}{|c|c|c|}
\hline
Grid Size & $\%$ time in fluid solver (RBF-IB) & $\%$ time in fluid solver (PL-IB) \\
\hline
{$64\times32$}   & 33.7  & 32.1 \\\hline 
{$128\times64$}   & 56.2  & 58.0     \\\hline 
{$256\times128$} & 65.4  & 79.3\\\hline
\end{tabular}	
\caption{Percentage of time per time-step spent in fluid solver as a function of grid size by both methods for $N_{\rm p} = 60$ platelets. The percentages for the RBF-IB method are the same for both $\Nd = 25$ and $\Nd = 50$ data sites, with the total time for the latter being larger. All results use $\Ns = 100$ sample sites (or IB points in the PL-IB method) per platelet.}	
\label{tab:Profile}
\end{table}

\subsection{Platelet Aggregation}

We now present the results of a platelet aggregation simulation with the RBF-IB simulation. We used the same boundary conditions, domain size, fluid properties and Poisseuille flow as in the previous subsection, but allow platelets to form links with other platelets and a portion of the chamber wall ($x = 0.4$ to $x = 0.7$) at the sample sites ($\Ns = 100$ per platelet). We used $\Nd = 50$ data sites per platelet, making the data sites a subset of the sample sites for convenience of visualization, and then visualize the data sites and the links between sample sites. We allowed each platelet to form up to $10$ links in total, either with the wall or with a neighbor; we allow links to cross each other for the purpose of simplicity, though this is usually prohibited in a platelet simulation. The simulation was run on a $128\times64$ grid with a time-step of $\Delta t = 10^{-4}$.

\begin{figure}[htbp]
\centering
\includegraphics[width=5.0in]{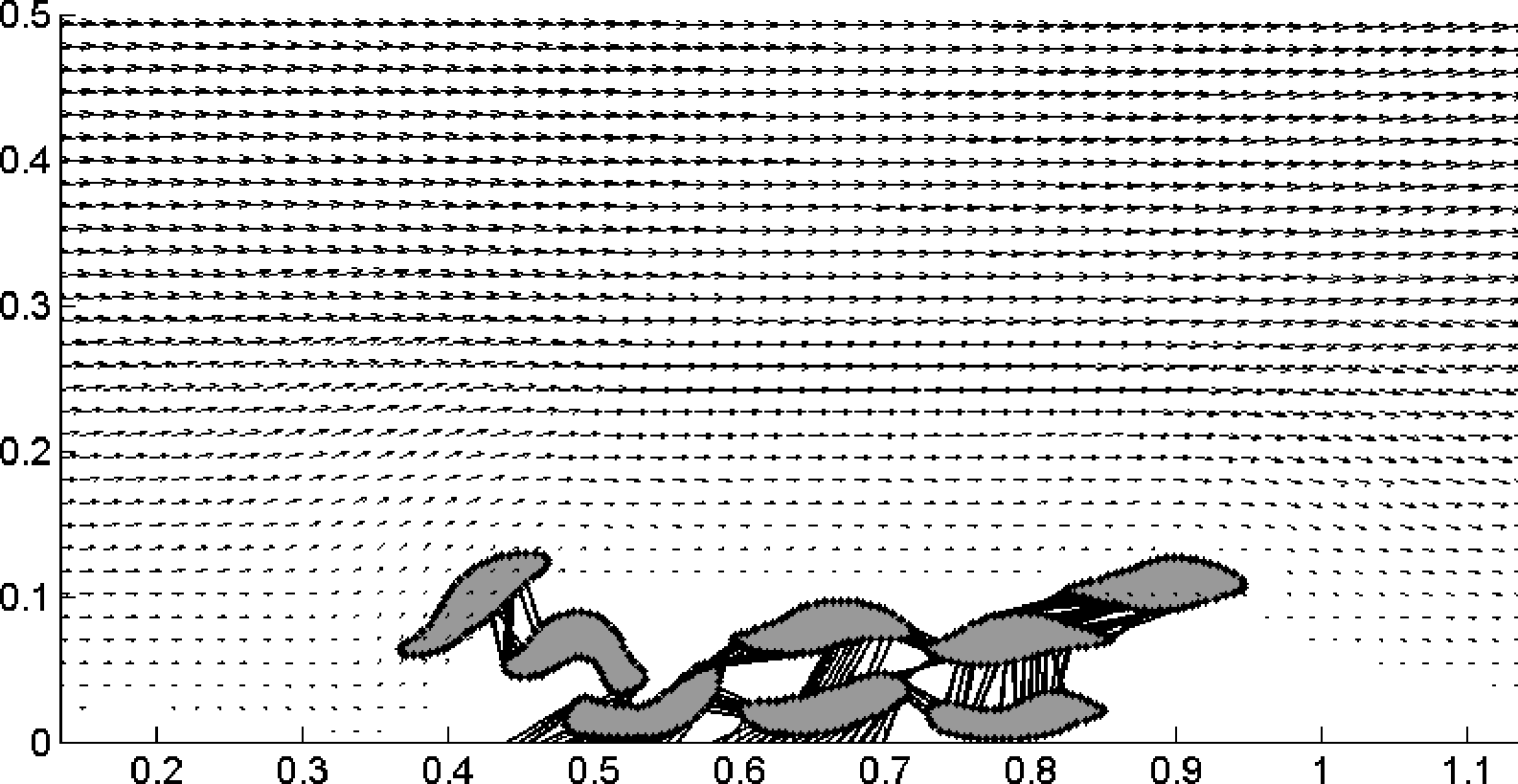}
\caption{Results of a platelet aggregation simulations with the RBF-IB method. The figure shows a zoomed-in snapshot of a platelet aggregation simulation achieved with the RBF-IB method with a time-step of $\Delta t = 10^{-4}$. The snapshot was taken at simulation time $t = 2.4$. The simulation was run on a $128\times64$ grid on a $[0,2]\times[0,1]$ domain. The arrows show the magnitude and direction of the velocity field.}
\label{fig:aggr1}
\end{figure}

Each platelet is initially an ellipse with radii $a = 0.06$ and $b = 0.015$. We initialize the platelets so that their centers of masses are at locations $(0.5,0.02)$, $(0.64,0.02)$, $(0.78,0.02)$, $(0.55, 0.07)$, $(0.68, 0.07)$, $(0.4, 0.045)$, $(0.23, 0.045)$ and $(0.65, 0.14)$. We chose these locations to ensure that three platelets lay on the wall, with three close enough to bind to the three bound to the wall, and two slightly further away. Each platelet attempts to maintain its initial elliptical shape. We then started the simulation and ran it to time $t = 2.4$. The results are shown in Figure \ref{fig:aggr1}. The figure shows both the velocity field and the platelet aggregate for a portion of $[0,2]\times[0,1]$ domain, the data sites on each platelet and the links between the sample sites corresponding to those particular data sites on the platelet.

There are two interesting features in Figure \ref{fig:aggr1}. The first is that the fluid flow gets diverted around the growing aggregate, a consequence of the size of the aggregate and the dynamics of the problem that mimics what one would hope to see in a realistic platelet aggregation simulation. The second feature is that some platelets are quite deformed, \emph{e.g.}, the platelet with center of mass approximately at $(0.5,0.02$), or its neighbor above and to its left. This is a consequence (and function) of the stiffness of each platelet, the shear rate of the flow and the number of links we allow each platelet to form. Higher platelet stiffness, lower shear rates and/or fewer (or weaker) links would lead to less deformation. The breaking model for inter-platelet and platelet-wall links can also affect the mechanics of aggregation. We note that our RBF model did not run into any instabilities in this simulation even when we ran it out to a time at which all the platelets (except the three closest to the wall) had left the domain.

\section{Summary}
\label{sec:discussion}

The IB method, as a numerical methodology for applications involving fluid structure interactions, naturally lends itself to our problem of interest: simulating platelet aggregation during blood clotting. In this application, platelets are modeled as immersed elastic structures whose shapes change dynamically in response to blood flow and chemistry. In previous work \cite{Shankar2012}, we discussed several geometric representations for platelets and compared them to the representation used within the traditional IB method. We concluded that an RBF geometric model for platelets would prove advantageous in several ways.

In this work, we explored the ramifications of using the RBF geometric model within the IB method, and compared the behavior of this new RBF-IB method against that of the traditional IB method. We discussed the issue of selecting an appropriate shape parameter for the RBF-IB method. We then presented a series of convergence studies for measuring errors and convergence both in the velocity field and in the representation of the immersed elastic structure. We went on to compare the computational costs incurred by both methods in the context of platelet simulations. We then compared the area conservation properties of both methods and also the time-step restrictions on both. We also remarked on the energy properties of our method.

We conclude the following:
\begin{itemize}
\item The RBF-IB method demonstrates first-order convergence, similar to that shown by the traditional IB method.
\item The RBF-IB method can be safely used both with a fixed number of data sites and with an increasing number of data sites as the background Cartesian grid is refined; the correct strategy is likely problem-specific;
\item Through the use of a global interpolant and the dual representation (data sites and sample sites), the RBF-IB method allows for a decoupling of accuracy in geometric modeling from accuracy in the full fluid solve, allowing coarse Lagrangian representations when sufficient.
\item The RBF-IB method is more computationally efficient than the traditional IB method, both due to the utilization of a small number of data sites and due to smoother forces being spread into the fluid resulting in a faster convergence from the fluid solver; and
\item The RBF-IB allows for larger time-step sizes than those allowed in the traditional IB method for a given grid size.
\end{itemize}

In previous work~\cite{Newren2007}, a sufficient condition for unconditional stability of an implicit IB method was established. The proof relied on the assumption that the set of points from which IB forces are spread is the same as that to which grid velocities are interpolated to update IB point positions. The RBF-IB method does not meet that condition, and it remains to be seen how this would impact an implict version of our method. Finally, an issue with the RBF-IB method is that it is dependent on the parametrization of the immersed elastic objects. For objects that are not easily parameterized in terms of circles and ellipses, the use of the RBF model as presented in our work (wherein the RBFs are restricted to the circle) may not be ideal. In the future, we thus hope to explore the use of RBFs in a meshfree variational form within the IB method so as to be able to easily evaluate constitutive models on arbitrary shapes. 

\vspace{0.2in}

{\bf Acknowledgments:}
We would like to acknowledge useful discussions concerning this work 
within the CLOT group at the University of Utah, with Professor Boyce Griffith
of New York University and with Professor Robert Guy of the University
of California, Davis.  The first, third and
fourth authors were funded under NIGMS grant R01-GM090203.  
The second author acknowledges funding support under NSF-DMS grant 
0540779 and NSF-DMS grant 0934581.

\appendix
\section{Algorithm for computing platelet forces with RBFs}
\label{sec:app1}

We now describe the implementation of the constitutive models outlined in Section \ref{sec:GeometricModeling}. We present algorithms for computing platelet forces in 2D.  

\underline{Notation}: In the description of the algorithms below we use standard matrix-vector operations such as multiplication as well as non-standard operations like element-by-element multiplication of matrices and vectors (sometimes called the Hadamard product).  We denote this latter operation with the $\circ$ operator.  For example, if $\vec{\boldsymbol{J}}$ and $\vec{\boldsymbol{L}}$  are $\Nd \times 2$ matrices and $\vec{R}$ is a vector of length $\Nd$ then the $i^{\text{th}}$ row of $\vec{\boldsymbol{J}}\circ\vec{\boldsymbol{L}}$ and $\vec{R}\circ\vec{\boldsymbol{J}}$ are given by
\begin{align*}
(\vec{\boldsymbol{J}}\circ\vec{\boldsymbol{L}})_{i,1:2} &= [(\vec{\boldsymbol{J}})_{i,1}(\vec{\boldsymbol{L}})_{i,1},\; (\vec{\boldsymbol{J}})_{i,2}(\vec{\boldsymbol{L}})_{i,2}] \\ 
(\vec{R}\circ\vec{\boldsymbol{L}})_{i,1:2} & = [(\vec{R})_{i}(\vec{\boldsymbol{L}})_{i,1},\; (\vec{R})_{i}(\vec{\boldsymbol{L}})_{i,2}]
\end{align*}
where $i=1,\ldots,\Nd$.

We define $\utngtd = \Dlond{1} \uvXd(t)$, the $\Nd \times 2$ matrix of tangent vectors at the data sites at time $t$ and $\vec{\|\tngtd\|}$, the $\Nd$ vector containing the two-norm of each row of $\utngtd$. The algorithm for computing platelet elasticity is as follows:

\begin{enumerate}
	\item Initialization ($t=t_0$): After creating and storing the RBF evaluation matrix as in Equation \eqref{eq:rbf_eval} and differentiation matrices as in Equations \eqref{eq:diff_mat_lond} and \eqref{eq:diff_mat_lons}, compute for each platelet:
		\begin{enumerate}
			\item The rest lengths for the tension force at the data sites: $\vec{l_0} = \utngtd = \Dlond{1} \uvXd(t_0)$. 
		  \item The bending-resistant force term for the platelet's initial configuration at the data sites, $\Dlons{4}\uvXd(t_0)$.
		\end{enumerate}
	\item For each time step ($t=t_k$, $k\geq 1$), compute for each platelet:
	 	\begin{enumerate}
	 	  \item The length of the tangent vectors $\utngtd = \Dlond{1} \uvXd(t_k)$ at the data sites: $\vec{\|\tngtd\|}$; and the unit tangents at the data sites: $\uutngtd$. 
	          \item The tension at the data sites, using the constitutive model: $\uTd = \kt(\vec{\|\tngtd\|} - \vec{l_0})$.  
		  \item The tension force at sample sites: $\uvFTs = \Dlons{1}\uvZd$, where $\uvZd = \uTd\circ\uutngtd$.
		  \item The bending force at sample sites: $\uvFBs = -\kb\lf(\Dlons{4}\uvXd(t_k) - \Dlons{4}\uvXd(t_0)\rt)$.
		  \item The interplatelet cohesion force from Equation \eqref{eq:force_IP} at the sample sites: $\uvFCs$. 
		  \item The total Lagrangian force at the sample sites: $\uvFs = \uvFTs + \uvFBs + \uvFCs$.
		\end{enumerate}
\end{enumerate}

\section{Time-stepping for the PL-IB method}

Here, we present the steps of the traditional Immersed Boundary algorithm when used with the RK2 time-stepping scheme from \cite{DevendranPeskin}.

\begin{enumerate}
\item Advance the structure to time level $t^{n + 1/2}$ using the current velocity field on the grid $\vu^n_g$. This is done by updating each IB point $\vX_q$ (for each q) using the equation
\begin{equation}
\vX^{n + 1/2}_q = \vX^n_q + \frac{\Delta t}{2} \vU^n_q \equiv \vX^n_q + \frac{\Delta t}{2} \sum_g \vu^{n}_g \delta_h (\vx_g - \vX^n_q) h^2,
\label{eq:IBmove1}
\end{equation}
where $h$ is the fluid grid spacing and $\delta_h$ is a discrete approximation to a
two-dimensional $\delta$-function. Here, $\vx_g$ and $\vX^{n+1/2}_q$ are the coordinates of grid point $g$ and IB
point $q$, respectively

\item The resultant $\vF^{n+1/2}_q$ of all of the force contributions that
act on an IB point $\vX^{n+1/2}_q$ is calculated for each $q$.  

\item
These forces are distributed to the Eulerian grid used for the
fluid dynamics equations using a discrete version of
Equation \eqref{eq:fluidf}:
\begin{equation}
\vf^{n+1/2}_g \equiv \vf^{n+1/2}(\vx_g) = \sum_q \vF^{n+1/2}_q \delta_h (\vx_g - \vX^{n+1/2}_q) dq.
\label{eq:spreadforce1}
\end{equation}
Here, $\vF^{n+1/2}_q$ is the Lagrangian force (per unit $q$) on
the IB point, $dq$ is the increment in parameter $q$ between consecutive
discrete sample sites, and $\delta_h$ is the same approximate $\delta$-function as used in Equation \ref{eq:IBmove1}.

\item
With the fluid force density $\vf^{n+1/2}_g$ now known at each grid point, the
fluid velocity is updated taking a half step ($\Delta t/2$) with a discrete
Navier-Stokes solver. As in \cite{DevendranPeskin}, we use a fractional-step projection method. First, a backward Euler discretization of the momentum equations is used. The pressure that enforces discrete incompressibility is determined~\cite{GUY:2005:SAP}. This gives us the velocity field $\vu^{n+1/2}_g$, the mid-step approximation required in an RK2 method.

\item
Using the mid-step fluid velocity $\vu^{n+1/2}_g$ and the mid-step IB point positions $\vX^{n+1/2}_q$, update the IB points $\vX^n_q$ for each $q$ to the time level $t^{n+1}$ using
\begin{equation}
\vX^{n+1}_q = \vX^n_q + \Delta t \vU^{n+1/2}_q \equiv \vX^n_q + \Delta t \sum_g \vu^{n+1/2}_g \delta_h (\vx_g - \vX^{n+1/2}_q) h^2,
\label{eq:IBmove2}
\end{equation}
where $\delta_h$ is the same approximate $\delta$-function we have used throughout.

\item
Update the velocity $\vu^n_g$ to time-level $t^{n+1}$ using the mid-step velocity $\vu^{n+1/2}_g$ and the force $\vf^{n+1/2}_g$. The mid-step velocities are advected, while a Crank-Nicolson scheme is used for time-stepping the momentum equations. The pressure projection gives us the discretely-incompressible velocity field $\vu^{n+1}_g$. Note that this step could have been performed as soon as $\vu^{n+1/2}_g$ was computed. It is independent of step (5).

\end{enumerate}

\section{Time-stepping for the RBF-IB method}

The RBF-IB method is time-stepped using the same RK2 method as above, with a few changes to incorporate data sites and sample sites.

\begin{enumerate}
\item Advance the structure to time level $t^{n + 1/2}$ using the current velocity field $\vu^n$. This is done by updating the \emph{data sites} $({\vX_{\rm d}})^n_j$ by a discrete analog of Equation \eqref{eq:ibU}
\begin{equation}
(\vX_{\rm d})^{n+1/2}_j = (\vX_{\rm d})^n_j + \frac{\Delta t}{2} (\vU_{\rm d})^n_j \equiv (\vX_{\rm d})^n_j + \frac{\Delta t}{2} \sum_g \vu^n_g
\delta_h (\vx_g - (\vX_{\rm d})^n_j) h^2.
\label{eq:IBmove3}
\end{equation}

\item Generate a new set of sample sites $\uvXs(t_{n+1/2})$ by applying the RBF evaluation operator to the data sites $\uvXd^{n+1/2}:=\uvXd(t_{n+1/2})$, {\em{i.e.}},
\begin{equation}
\uvXs(t_{n+1/2}) =  \Es\uvXd(t_{\rm new}).
\end{equation}

\item The total force at the sample sites $\uvFs^{n+1/2}$ is calculated using the algorithm from Appendix A.

\item
These forces are distributed to the Eulerian grid used for the
fluid dynamics equations using a discrete version of
Equation \eqref{eq:fluidf}:
\begin{equation}
\vf^{n+1/2}_g \equiv \vf^{n+1/2}(\vx_g) = \sum_q \vF^{n+1/2}_q \delta_h (\vx_g - (\vX_{\rm s})^{n+1/2}_q) dq.
\label{eq:spreadforce2}
\end{equation}
Here, $\vx_g$ and $(\vX_{\rm s})^{n+1/2}_q$ are the coordinates of grid point $g$ and sample site $q$, respectively, $\vF^{n+1/2}_q$ is the Lagrangian force (per unit $q$) on the sample site, $dq$ is the increment in parameter $q$ between consecutive
discrete sample sites, and $\delta_h$ is the same approximate $\delta$-function as used in Equation \ref{eq:IBmove3}.

\item
With the fluid force density $\vf^{n+1/2}_g$ now known at each grid point, the
fluid velocity is updated taking a half step ($\Delta t/2$) with a discrete
Navier-Stokes solver. Again, we use a fractional-step projection method, with a backward Euler discretization of the momentum equation, and a projection to determine the pressure that enforce incompressibility~\cite{GUY:2005:SAP}. This gives us the velocity field $\vu^{n+1/2}_g$, the mid-step approximation required in an RK2 method.

\item
Using the mid-step fluid velocity $\vu^{n+1/2}_g$ and the mid-step \emph{data site} positions $(\vX_{\rm d})^{n+1/2}_j$, update the \emph{data sites} $(\vX_{\rm d})^{n+1/2}_j$ for each $j$ to the time level $t^{n+1}$ using
\begin{equation}
(\vX_{\rm d})^{n+1}_j = (\vX_{\rm d})^n_j + \Delta t (\vU_{\rm d})^{n+1/2}_j \equiv (\vX_{\rm d})^n_j + \Delta t \sum_g \vu^{n+1/2}_g \delta_h (\vx_g - (\vX_{\rm d})^{n+1/2}_j) h^2,
\label{eq:IBmove4}
\end{equation}
where $\delta_h$ is the same approximate $\delta$-function we have used throughout.

\item
Update the velocity $\vu^n_g$ to time-level $t^{n+1}$ using the mid-step velocity $\vu^{n+1/2}_g$ and the force $\vf^{n+1/2}_g$. The mid-step velocities are advected, while a Crank-Nicolson scheme is used for time-stepping the momentum equations. The pressure projection gives us the discretely-incompressible velocity field $\vu^{n+1}_g$. 

\end{enumerate}

Observe that the data sites are updated twice per time-step in the RK2 scheme, but the sample sites are only generated once. Since the data sites are typically a fraction of the number of IB points from the PL-IB method, the computational cost is significantly lower for the RBF-IB method, even factoring in the interpolation and the sample site generation.

\bibliographystyle{wileyj}
\bibliography{flddoc}

\end{document}